\documentclass[paper=A4,10pt]{article}

\usepackage{algorithm}
\usepackage{algorithmic,eqparbox,array}
\usepackage{amsmath}
\usepackage{amsthm}
\usepackage{amssymb}
\usepackage{natbib}
\usepackage{geometry}
\usepackage[adobe-utopia]{mathdesign}
\usepackage[T1]{fontenc}
\usepackage{pdflscape}
\usepackage{booktabs}
\usepackage{hyperref}
\usepackage[textwidth=2.5cm]{todonotes}
\usepackage{authblk}
\usepackage{pbox}
\usepackage{pifont}

\hypersetup{colorlinks, urlcolor=green, linkcolor=red,
citecolor=blue, bookmarksnumbered}

\newcommand{\converge}{\pbox{1em}{$\downarrow$\\ $\uparrow$}}
\newcommand{\diverge}{\pbox{1em}{$\uparrow$\\ $\downarrow$}}

\newcommand{\eqn}[1]{\begin{equation} #1 \end{equation}}

\newcommand{\R}{\mathbb{R}}

\newcommand{\constant}{\ensuremath{\text{c}}}

\newcommand{\rv}{\varepsilon}

\newcommand{\prob}{\text{P}}

\newcommand{\scale}{\mu}

\newcommand{\Kmin}[2]{\LL_{#1\setminus#2}}

\newcommand{\rs}{R}

\newcommand{\LL}{\mathcal{L}}
\newcommand{\y}{\boldsymbol{y}}
\newcommand{\z}{\boldsymbol{z}}
\newcommand{\RV}{\boldsymbol{\varepsilon}}
\newcommand{\ps}{\ensuremath{\text{PS}}}
\newcommand{\cf}{\ensuremath{\varphi}}
\newcommand{\ic}{\ensuremath{\alpha}}
\newcommand{\MD}{\ensuremath{{\text{M}\Delta\text{-}}}}
\newcommand{\Md}{\ensuremath{{\text{M}\Delta}}}
\newcommand{\M}{\ensuremath{\text{M-}}}
\newcommand{\Mm}{\ensuremath{\text{M}}}
\newcommand{\A}{\ensuremath{\text{A-}}}
\newcommand{\Aa}{\ensuremath{\text{A}}}
\newcommand{\MEV}{\ensuremath{\text{MEV}}}

\newcommand{\GMEV}{\ensuremath{\text{GMEV}}}
\newcommand{\MNL}{\ensuremath{\text{MN}}}
\newcommand{\PSL}{\ensuremath{\text{PS}}}
\newcommand{\PCL}{\ensuremath{\text{PC}}}
\newcommand{\LNL}{\ensuremath{\text{LN}}}

\newcommand{\pspar}{\beta}
\newcommand{\yA}{\ensuremath{\y^\text{A}}}
\newcommand{\yM}{\ensuremath{\y^\text{M}}}

\newcommand{\yMDr}{\ensuremath{\y^{\text{M}\Delta,r}}}
\newcommand{\yyMD}{\ensuremath{y^{\text{M}\Delta}}}
\newcommand{\yyMDr}{\ensuremath{y^{\text{M}\Delta,r}}}
\newcommand{\Pref}[1]{P^\text{ref}(#1)}
\newcommand{\GMNL}{\ensuremath{G^\text{MN}}}
\newcommand{\GPSL}{\ensuremath{G^\text{PS}}}
\newcommand{\GPCL}{\ensuremath{G^\text{PC}}}
\newcommand{\GLNL}{\ensuremath{G^\text{LN}}}
\newcommand{\one}{\boldsymbol{1}}
\newcommand{\TT}{{\ensuremath{\tau}}}
\newcommand{\demand}{D}
\newcommand{\flow}{f}
\newcommand{\flowvec}{\boldsymbol{f}}
\newcommand{\flowfeas}{\varOmega}
\newcommand{\kkta}{\kappa}
\newcommand{\kktb}{\lambda}
\newcommand{\cmark}{\ding{51}}%
\newcommand{\xmark}{\ding{55}}%

\DeclareMathOperator{\E}{\mathbb{E}}
\DeclareMathOperator{\Var}{Var}
\DeclareMathOperator{\stdev}{\sigma}
\DeclareMathOperator{\Covar}{Cov}

\begin{document}

%
% TITLE PAGE, including submission date, word count (max 7500 'words', floats (figures and tables) each count as 250 words), and author names, affiliations, addresses, phone numbers, fax numbers, and e-mails (please indicate corresponding author)
%

\title{Generalized Multivariate Extreme Value Models for Explicit Route Choice Sets}
\author[1,2*]{Erik-Sander Smits}
\author[1]{Adam J. Pel}
\author[3]{Michiel C.J. Bliemer}
\author[1]{Bart van Arem}

\affil[1]{\small Department of Transport \& Planning, Delft University of Technology, Postbus 5048, 2600 GA Delft, The Netherlands}
\affil[2]{\small Arane Adviseurs in Verkeer en Vervoer, Groen van Prinsterersingel 43b, 2805 TD Gouda, The Netherlands}
\affil[3]{\small Institute of Transport and Logistics Studies, The University of Sydney, NSW 2006, Australia}
\affil[*]{\small Corresponding author: e.smits@arane.nl, +31 (0)6 10857379}

\maketitle

\begin{abstract}  
This paper analyses a class of route choice models with closed-form probability expressions, namely, Generalized Multivariate Extreme Value (GMEV) models. A large group of these models emerge from different utility formulas that combine systematic utility and random error terms. Twelve models are captured in a single discrete choice framework. The additive utility formula leads to the known logit family, being multinomial, path-size, paired combinatorial and link-nested. For the multiplicative formulation only the multinomial and path-size weibit models have been identified; this study also identifies the paired combinatorial and link-nested variations, and generalizes the path-size variant. Furthermore, a new traveller's decision rule based on the multiplicative utility formula with a reference route is presented. Here the traveller chooses exclusively based on the differences between routes. This leads to four new GMEV models. We assess the models qualitatively based on a generic structure of route utility with random foreseen travel times, for which we empirically identify that the variance of utility should be different from thus far assumed for multinomial probit and logit-kernel models. The expected travellers' behaviour and model-behaviour under simple network changes are analysed.  Furthermore, all models are estimated and validated on an illustrative network example with long distance and short distance origin-destination pairs. The new multiplicative models based on differences outperform the additive models in both tests.       
\end{abstract}


\section{Introduction}

Route choice is important in transport applications such as network equilibrium modelling, day-to-day route choice decisions, and route guidance. The literature describes various methods to model the route choice behaviour of travellers. They range from simple deterministic shortest route choice to sophisticated stochastic models in which random error terms capture travel time uncertainty and taste heterogeneity amongst travellers. Commonly, the route choice model is applied in an iterative process, for example, to reach equilibrium in a congested network or for en-route decisions. For large networks the number of times route choice has to be simulated is very high. Therefore, a good balance between realism and computational efficiency is required.

\subsection{Random utility maximization}

A common choice mechanism is that travellers consider multiple routes, assign a subjective utility to each route, and choose the route with the highest utility. The discrete choice framework based on random utility maximization (RUM) allows random components in the utility formulations. In this case the utility of each route follows some random distribution. Observe that a route consists of a set of links, and we assume that a route's utility is (among other things) determined by the characteristics of the corresponding links. Then the joint probability distribution of route utilities contains dependencies when routes have overlap, and thus share (dis)utility from the same link(s). Furthermore the distribution is heteroscedastic (i.e., the variability of utility differs amongst routes) since routes have different lengths. Section \ref{sec:rum} provides an in-depth analysis of desired random route utility properties, that especially explores the distribution of foreseen travel times besides the usual analyst error. This gives new insights into the variance and covariance structure of route utility.

The Multinomial Probit (MNP) model for route choice \citep[see ][]{Daganzo01081977,Yai1997} can easily address correlation and heteroscedasticity, but no closed-form formulation of the route probabilities exists. The latter leads to computationally expensive simulations. On the other hand the Multinomial Logit (MNL) model for route choice \citep[see ][]{Daganzo01081977}\footnote{\citet[][Eqn. 10]{Daganzo01081977} are the first to present MNL route choice probabilities, based on the method of \citet{Dial197183}} has an elegant closed-form formulation for the route probabilities; however, route utilities have to be independent and homoscedastic. Extensions to mixed logit models \citep{Bekhor2002} are possible to overcome some of the limitations, however, these models again have to rely on numerical solution methods requiring simulations which are often infeasible on large scale networks. In addition, the thus far proposed covariance structures in MNP and mixed logit are not in line with revealed foreseen travel times; the \emph{variance} has been assumed proportional to the mean, whereas instead, data shows a linear relation between \emph{standard deviation} and mean (see Section \ref{sec:rum:lin}). 

Several adaptations to the MNL model exist to address correlation due to route overlap \citep{Prato2009}. They either insert a correction term into the utility, or exploit the more general Multivariate Extreme Value (MEV) distribution for the error term. In Section \ref{sec:instances} an overview of these methods is presented. The heteroscedasticity is less addressed in the literature. In practice, the scale parameter of MNL is sometimes considered to be proportional to the distance between origin and destination, but this does not address heteroscedasticity within the origin-destination (OD) route set \citep[see ][]{Chen2012}. Less pragmatic is the solution by  \citet{Castillo2008373} who assume Weibull distributed utilities leading to the Multinomial Weibit (MNW) model. Only slightly different is the approach of \citet{Galvez2002,Fosgerau2009494} who use a multiplicative error term instead of an additive error term. The probabilities of the earlier-known Kirchhoff distribution of routes as presented by \citet{Fellendorf2010} also coincide with MNW. Recently, \citet{Kitthamkesorn2013} included the path-size factor into the MNW model to correct for correlation between routes. 

\subsection{Generalized Multivariate Extreme Value models}

In this paper the state-of-the-art of models with closed-form utility formulations is reviewed, and they are gathered in a framework with a single, general form for the choice probabilities. Two classes of MEV can be identified -- each based on a different utility formula --, which are unified in what we call Generalized Multivariate Extreme Value (GMEV). The generalization gives rise to three undiscovered route choice models based on the multiplicative utility formula. Furthermore, we recall that equivalent MEV formulations exist for models in which a correction term is inserted in the utility formulation (e.g. Path-Size Logit). All these models can therefore be captured in the GMEV framework based on either additive or multiplicative utility formulas. Section  \ref{sec:qualassess} assesses these models qualitatively by analysing their distributions and comparing them with the desired properties formulated in Section \ref{sec:rum:req}. Additive models prove not to be able to capture random foreseen travel time, but they do capture the analyst error. On the other hand, multiplicative models can handle both errors, but based on the same distribution (i.e., the foreseen travel time and analyst error are completely dependent). 

One of the desired properties cannot be fulfilled by any existing model. Namely, in a transport network model the route choice model is applied to many route sets, so it is important that the model gives realistic results for all route sets existing in a network. In particular, the results should remain realistic when the network and/or routes change. With respect to this, both MNL and MNW have an undesired property. In MNL, the route probabilities will remain equal if a constant is added to each route's utility. In MNW, the route probabilities will remain equal if each route's utility is multiplied with a constant. Both these properties are unrealistic. Therefore, \citet{Xu2015} provide a hybrid method where the choice probabilities have a MNL en MNW component; however, this hybrid method stems from a combined impedance function (i.e., it defines the choice probabilities based on the choice probability functions of logit and weibit). Although the model can be useful in practice, the model does not stem from a utility formulation with an additive and multiplicative error term (Section \ref{sec:hybrid} discusses how the hybrid models would fit in the \GMEV\ framework). 

To partially overcome the shortcomings of MNW, we introduce new \GMEV\ route choice models based on the multiplicative utility formula and a reference route. The probability of choosing the reference route only depends on the non-overlapping differences with other routes, since the overlap is explicitly removed. The models fit in the same framework as the existing models. Furthermore, they are designed to be suitable under changing networks and routes.   A qualitative assessment and a numerical benchmark on a small network example of the existing and new models show that the new models can resemble expected behaviour under more types of network changes than existing closed-form models can. For example, the model performs well if all route costs are multiplied with a factor, but also if a constant cost is added to all routes. 

\subsection{Route set generation}

This paper assumes that a relevant route set is explicitly given. For an overview of route set generation techniques, see \citet{Frejinger2009984,Prato2012}. However, it should be noted that the generation and composition of the route set has a large influence on model outcomes \citep{bliemertrr2008,Cascetta2002,Prato2012}.  Another disadvantage with explicit route sets is that a correct sample of routes is required to obtain unbiased parameter estimates \citep{Frejinger2009984}.  For discussions on sampling routes based on distributions derived from route choice models see \citep{Frejinger2009984,Flotterod2013,Guevara2013}. \citet{Bierlaire2008381} points out difficulties with selection bias in the estimation of additive MEV models from choice based samples.  A recently revisited different approach -- that overcomes route set generation/sampling related problems -- is to implicitly generate routes as done by \citet{Dial197183,Papola2013CACIE,Fosgerau201370,Mai2015}. However, these models either have a restricted route set \citep{Dial197183,Papola2013CACIE} (see Section \ref{sec:jointnetwork}), or contain all routes, even those with loops \citep{Fosgerau201370,Mai2015}. In addition, no methods with implicit choice set generation that can handle non-additive travel costs are known to exist by the authors of this paper. 

Since explicit finite route sets have practical advantages, they are commonly used in traffic assignment models to either avoid too computationally expensive (dynamic) shortest-path calculations, to allow dynamic network loading of path flows, or to allow for traffic assignment models with non-additive link costs. For example, \citet{Zhou2015} generates an a-priori route set in a static assignment context using the deterministic method of \citet{Bekhor2006}, \citet{tr-b14} generates an a-priori route set for a quasi-dynamic model using the stochastic method of \citet{Fiorenzo2004}.    

\subsection{Contribution}

The contribution of this paper is summarized as follows. First, we describe the desired properties of random route utility which includes random foreseen travel times, and empirically verify the linear relation between mean and standard deviation of foreseen travel time. Second, existing closed-form route choice models based on random utility maximization are presented in a novel \GMEV\ framework and it is shown that all existing closed-form models with explicit choice sets fit in the \GMEV\ framework. Third, the unexplored area of multiplicative MEV models gives rise to three new closed-form route choice models, and the proposed multiplicative MEV models based on a reference route give rise to another four new closed-form route choice models. The working of this model is illustrated with a small network example, and its behaviour with respect to network changes  is compared to the behaviour of existing models. Fourth, we show that the additive models, contrary to multiplicative models, can not capture random foreseen travel time. Fifth, we show the strength of a unified GMEV model by providing its stochastic user equilibrium formulation. Sixth, all twelve models (of which seven are new) are estimated on a carefully constructed small network example with multiple choice sets for origin-destination pairs with both short and long distances. The estimation is done twice with two different synthetic datasets generated by MNP simulations based on the desired properties of random route utility, this allows benchmarking the models on the other dataset. The benchmark provides insight in the performance of the model if its applied on different networks, without re-estimation. Overall, when the models are estimated on one synthetic dataset and validated on the other synthetic dataset, the new models based on reference routes have a better validation than existing models. They can better approximate -- without simulations -- the MNP probabilities that take, amongst others, heteroscedasticity and correlation into account.

\section{Random Route Utility Formulation}\label{sec:rum}

The random utility maximisation (RUM) framework is commonly used for describing route choice. Each route alternative is associated with some utility and the traveller is assumed to choose the route that provides the maximum expected utility. Since utility is not directly observable, analysts typically assume a structure of utility that includes measurable components and random terms that describe unmeasurable components. This section provides a route utility formulation with random foreseen travel times, that is used in Section \ref{sec:qualassess} to assess the GMEV models. 

Let $\rs$ denote the set of relevant routes between a certain origin-destination (OD) pair. The systematic utility $V_r$ for a route $r\in\rs$ is decomposed in a part based on travel time $\TT_r$, and `other' systematic utility $V^0_r$, thus $V_r=V^0_r+\beta\TT_r$, where $\TT_r$ is a random variate representing foreseen travel time, 
%\footnote{Expected travel time should not be confused with the expected value of travel time. The first is a random distribution while the latter is the mean of that distribution. Also, the latter has nothing to do with expected value of time (VOT). The terminology expected travel time is used since when travellers make their route choice decision, they do that based on the travel times they expect for each route.}
$\beta$ is the travel time parameter. Then we write the random utility for route $r\in\rs$ as\footnote{The OD pair and individual are not indexed since they are not relevant for this study and omitting them makes the formulas more readable.}, 
\begin{equation}
\label{eqn:U}
U_r=V_r+\rv_r=V^0_r+\beta\TT_r+\rv_r, 
\end{equation}
where $\rv_r$ is the analyst error. $V^0_r$ is assumed to be deterministic and can contain all route attributes other than the foreseen travel time, such as travel distance, running cost, number of traffic lights encountered, number of left turns, travel time variability, etc. Typical random component $\rv_r$ includes attributes that are considered by travellers when choosing their routes, but that are not included in $V^0_r+\beta\TT_r$. Since the analyst does not have perfect knowledge on how the traveller gathers information about its travel time, we explicitly consider foreseen travel time as a random distribution $\TT_r$ in this work. In addition, travellers can have different preferences with respect to travel time and they do not perfectly measure the actual values of the attribute levels, but rather have a subjective perception that deviates from the actual values. Denote the analyst's estimate of foreseen travel time with $\hat\TT_r$, which generally does not consider all conditions (e.g., exact departure time, weather, information sources) the traveller is aware of. The deviation between the actual and estimate of foreseen travel time is typically larger for longer trips (see next Section). For example, when asking a traveller about the travel time of a recent short trip, his or her expectation may be off by one or two minutes, while for a long trip the expectation may be off by 10 minutes or more. We would like to point out that the foreseen travel time distribution does not capture travellers' disproportionate resilience towards high travel time variability, which can be captured by an (un)reliability term in $V^0_r$.

The utility in Equation (\ref{eqn:U}) is very general since the analyst can include any attribute in $V^0_r$, and the, in our opinion, two most important sources of randomness for route utility are captured. Extension for errors around other terms in the systematic utility can be made in a  straightforward fashion. We focus in this paper on the foreseen travel time because it is the most important factor in route choice. In addition, data is available to support assumptions on its structure.

Both $\TT_r$ and $\rv_r$ are random variables and are described by probability distributions. A common assumption for analyst errors $\rv_r$ is that they are independently and identically distributed (i.i.d.). The next section discusses the foreseen travel time distribution in more detail. Since $\TT_r$ and $\rv_r$ capture different sources of randomness, they can be considered independent, i.e. $\Covar(\TT_r,\rv_r)=0$ for all $r\in R$. 

\subsection{Structure of foreseen travel time}\label{sec:rum:lin}

Consider two similar road segments with lengths respectively 1 and 5 kilometres, and suppose the distribution of foreseen travel times for the 1-km road segment is known. There are two natural ways to determine the foreseen travel time distributions of the 5-km road segments by, which impose different relations between their means and variances. The first method -- used in all previous Probit studies \citep[e.g., ][]{Yai1997,Daganzo01081977} -- is based on link-additivity and postulates that the foreseen travel time of the 5-km road segment is equally distributed as the sum of foreseen travel times of five independent 1-km road segments. This implies for each set of road segments for which the \emph{additivity postulate} holds, that the ratio between the mean and the \emph{variance} is equal for each segment. Link-additivity is convenient in transport networks, since overlap and heteroscedasticity of routes can automatically be captured when the normal distributed travel costs are drawn for each link and then summed per route. The second method is based on scaling and postulates that the foreseen travel time of the 5-km road segment is equally distributed as five times the foreseen travel time of one 1-km road segment. This implies for each set of road segments for which the \emph{scaling postulate} holds, that the ratio between the mean and the \emph{standard deviation} is equal for each segment. Furthermore, the convenience of link-additivity does not apply under the scaling postulate. 

When an analyst wants to specify means and variances of foreseen travel times of routes, these two postulates each imply a different structure. Under the natural assumption that the mean equals the analyst's estimate $\hat\TT_r$, the variance is linear in $\hat\TT_r$ under the linear-additivity postulate and quadratic in $\hat\TT_r$ under the scaling postulate. In order to examine which postulate is more plausible, we first analyse the relationship between mean travel times and standard deviation of travel times in two datasets, and then discuss findings in the literature related to the postulates.

%The distribution of perceived travel time follows directly from the distribution of the perception error. Although we do not know the exact distribution of $\rvper_r$, without loss of generality $\E(\rvper_r) = 1$.\footnote{Any structural bias towards over- and underestimation can equivalently be captured in $\beta$} We assume that the $\rvper_r$'s are identically distributed -- but not necessary independent -- with $\Var(\rvper_r)=\theta^2$, where $\theta$ is the so called proportionality parameter, such that the standard deviation is $\stdev(\rvper_r)=\theta$. Then it follows that $\Var(\TT_r)=\hat \TT_r^2\theta^2$ and $\stdev(\TT_r)=\hat \TT_r\theta$; next, this linear relation between the estimate and standard deviation of travel time is underpinned with empirical data. 

\begin{figure}
\includegraphics[width=\textwidth]{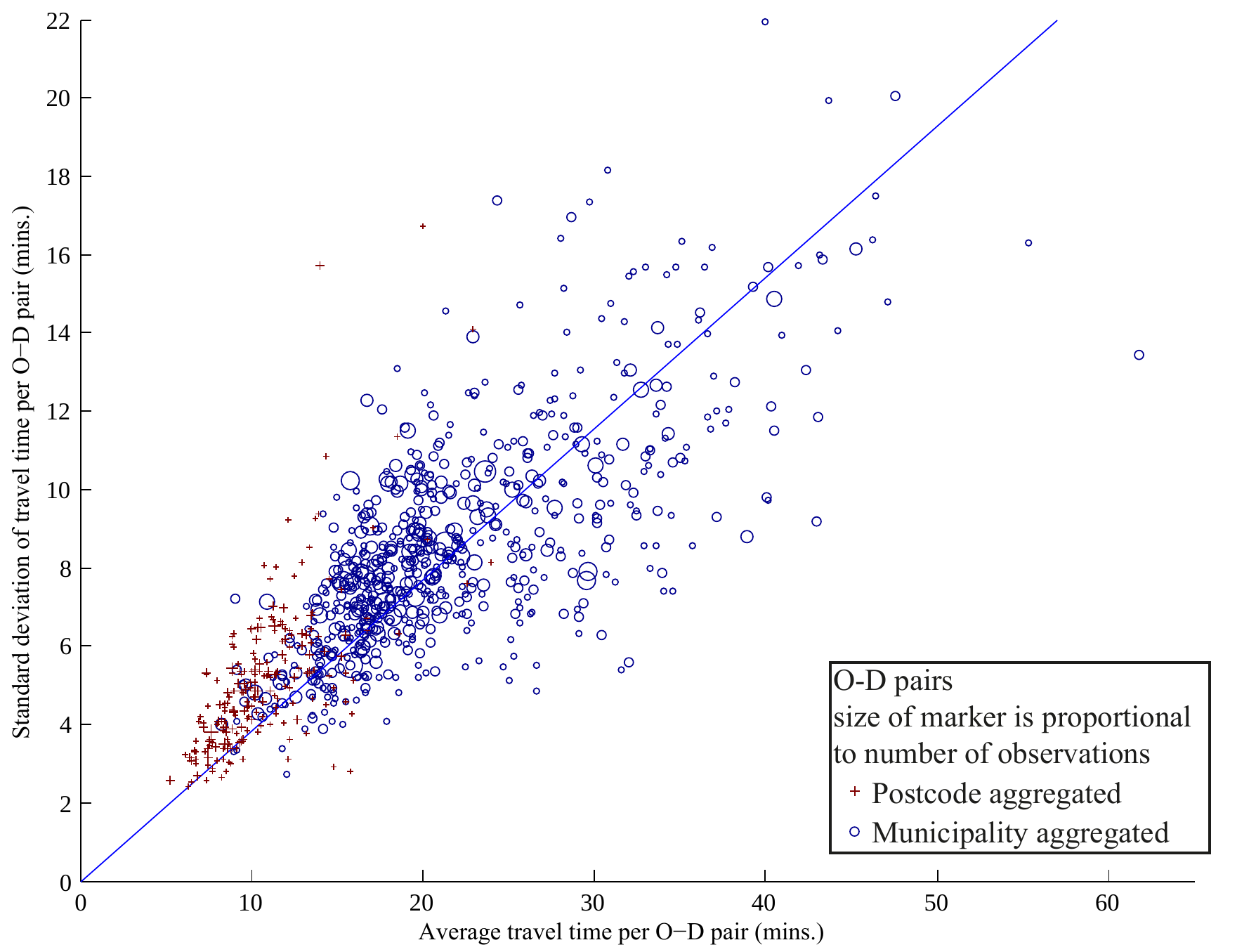}
\caption{Relation between the average and the standard deviation of revealed perceived travel times for an OD pair, shown for 216 postcode OD pairs and 686 municipality OD pairs from Dutch national travel survey data. Linear regression between mean and standard deviation has a R-square of 0.9379. Linear regression between mean and variance has a R-square of 0.8016. The linear regression line has slope 0.3859.}
\label{fig:TTscatter}
\end{figure}

\begin{figure}
\begin{center}
\includegraphics{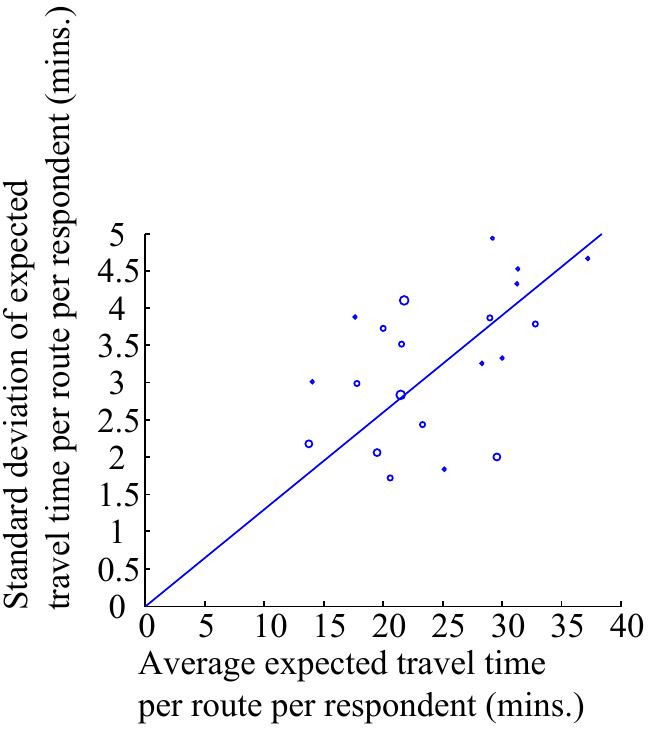}
\end{center}
\caption{Relation between the average and the standard deviation of foreseen travel times for a route from a respondent, shown for 21 routes from 19 respondents from trips between Delft and The Hague. The marker size is proportional to the number of entries per route. Linear regression between mean and standard deviation has an R-square of 0.9300. Linear regression between mean and variance has an R-square of 0.8401. The linear regression line has slope 0.1301.}
\label{fig:TTscatterSurvey}
\end{figure}

We have analysed household travel survey data from the Netherlands for the years 2010, 2011, and 2012 (Onderzoek Verplaatsingen in Nederland, OViN). In this survey, travellers are asked to state the travel times from trips they have made during a single day. We pooled data from trips with motorised private vehicles (i.e., car-driver, car-passenger, and motor bicycle) per 4-digit postal code and municipality\footnote{The four large municipalities Amsterdam, Utrecht, The Hague and Rotterdam are subdivided in respectively 7, 3, 5 and 5 districts.}. After excluding OD pairs with less than 25 observations we obtained 686 municipality aggregated OD pairs and 216 postal code aggregated OD pairs containing 38,106 and 8,079 trips, respectively. Figure \ref{fig:TTscatter} shows the relationship between the average travel time and the standard deviation of the perceived travel times for each OD pair. The data shows a strong correlation and a linear relation between the average travel time and the standard deviation of travel time. Under the assumption that the traveller had complete information, the perceived travel time and foreseen travel time coincide, and the data is then in accordance with the scaling postulate. Clearly, the spread in reported travel times here is also the result of variations in traffic conditions and travel demand, routing, and aggregation of different households in each OD pair. Unfortunately, we do not have the exact data on the foreseen travel time per route in the household travel surveys; however, we do have data from another smaller study about route choice and information.

This second data set comes from a study where 32 commuters between The Hague and Delft in the Netherlands were asked about the used information sources \citep[see ][]{Ramos2012revealed,Ramos2015thesis}. Some of them were equipped with GPS devices for real-time information. After every trip they made, they were asked to indicate what their expected route and foreseen travel time for their next trip are. We aggregated this data per route per respondent, and determined the mean and standard deviation for all routes with at least 10 entries. In total 407 entries for 21 routes from 19 respondents\footnote{Two respondents reported two routes at least 10 times; 17 other respondents reported one route at least 10 times.} are used for the results in Figure \ref{fig:TTscatterSurvey}. This dataset also indicates a linear relationship between the average and standard deviation. Thus both datasets indicate that the scaling postulate is more plausible than the additivity postulate.

\citet{MahmassaniTRR12} investigate travel time variability based on GPS and simulated trajectories and they conclude that a linear relation exists between the standard deviation of observed travel times and the average travel times. They furthermore reject a square root and quadratic relation. \citet{MahmassaniTRR13} provide the same conclusion at a network level, rather than just at link level. \citet{fosgerau2010relation} and \citet{yildirimoglu2013exploring} use data on time-of-day dependent relations between mean and standard deviation to show hysteresis. Their plots also strongly suggest a linear relation. Furthermore, \citet{hellinga2012} took a year of measurements from loop detector data on a motorway and aggregated each day into 15-minute averages. They also show a strong linear relationship between the standard deviation of the observed travel times and the average travel times. Despite that none of these studies analysed \emph{foreseen} travel time distributions, they all indicate that the scaling postulate holds for their travel times. These results, the OViN data, and the actual foreseen travel time data, should convince that the scaling postulate is more plausible than the additivity postulate for distributions of foreseen route travel times.

Write $\stdev(\TT_r)$ as the standard deviation of foreseen travel time. Let this standard deviation be proportional to $\hat \TT_r$ with proportionality constant $\theta$, then $\stdev(\TT_r)=\theta\hat{\TT_r}$  

This finding contradicts the assumptions for MNP route choice models made by \citet{Daganzo01081977} and \citet{Yai1997}. Their assumption of a linear relationship between the variance of travel times and $\hat \TT_r$, stems from the additivity postulate. In these   models it is practical that route travel time can be readily obtained from the summation of normal distributed link travel times, but this seems to be unrealistic. It is possible to adjust the MNP-model to make it compatible with the scaling postulate by directly determining the (co)variance-matrix for all routes. However, this comes at the cost of link-additivity.  

\subsection{Covariances of foreseen travel times}

Since routes may overlap, foreseen travel times are correlated. Route overlap induces a structure of covariances between travel times. This means that if two routes are largely overlapping, then also the foreseen travel time of both routes should be almost identical. Consider two routes, $r,s \in \rs$. Let $\hat \TT_r=\hat\TT_{rs}+\hat\TT_{r\setminus s}$ be the decomposed travel time estimates of route $r$, where $\hat\TT_{rs}$ is the overlapping part with route s, and $\hat\TT_{r\setminus s}$ the non-overlapping part. The covariance between $\TT_r$ and $\TT_s$ should capture the overlap. We propose two candidate formulations for the covariance. One based on a correlation formulation derived from the literature, and one based on variance distribution over links.

The approaches from \citet{Daganzo01081977,Yai1997,Bekhor2002} rely on the additivity postulate, and the additivity over links leads to a definition of covariances. Under the scaling postulate these covariances are not directly usable due to the different relation between mean and standard deviation. However, the corresponding correlation, $\hat\TT_{rs}/\sqrt{\hat\TT_r\hat\TT_s}$, is independent of any relationship between mean and standard deviation, and can be used together with $\stdev(\TT_r)=\theta\hat\TT_r$. Then from the definition of correlation it follows that
\begin{align}
\text{Corr}(\TT_r,\TT_s)&=\frac{\hat\TT_{rs}}{\sqrt{\hat\TT_r\hat\TT_s}}=\frac{\Covar(\TT_r,\TT_s)}{\stdev(\TT_r)\stdev(\TT_s)}\\
&\Downarrow\nonumber\\
\Covar(\TT_r,\TT_s)&=\theta^2\hat\TT_r\hat\TT_s\frac{\hat\TT_{rs}}{\sqrt{\hat\TT_r\hat\TT_s}}=\theta^2\hat\TT_{rs}\sqrt{\hat\TT_r\hat\TT_s}.\label{eqn:covGm}
\end{align}

An alternative formulation of the covariance can be obtained by using the definition $\Var(\TT_r-\TT_s)=\Var(\TT_r)+\Var(\TT_s)-2\Covar(\TT_r,\TT_s)$ and an approximation of variance of the non-overlapping parts of travel time. Under the assumption that a route's variance is equally `distributed' over the total travel time and that $\TT_r$ can be decomposed in independent parts, $\TT_r=\TT_{rs}+\TT_{r\setminus s}$, the variance associated with the non-overlapping part equals 
\begin{equation}
\Var(\TT_{r\setminus s})=\frac{\hat\TT_{r\setminus s}}{\hat\TT_r}\Var(\TT_r)=\frac{\hat\TT_{r\setminus s}}{\hat\TT_r}\theta^2\hat\TT_r^2=\theta^2\hat\TT_r\hat\TT_{r\setminus s}.
\end{equation}
Then the covariance equals
\begin{align}
	\Covar(\TT_r,\TT_s)&=\frac{1}{2}\left(\Var(\TT_r)+\Var(\TT_s)-\Var(\TT_r-\TT_s)\right)\\
	&=\frac{1}{2}\left(\Var(\TT_r)+\Var(\TT_s)-\Var(\TT_{r\setminus s}-\TT_{s\setminus r})\right)\\
	&=\frac{1}{2}\theta^2\left(\hat\TT_r^2+\hat\TT_s^2-\hat\TT_r\hat\TT_{r\setminus s}-\hat\TT_s\hat\TT_{s\setminus r}\right)\\
	&=\frac{1}{2}\theta^2\left(\hat\TT_r(\hat\TT_r-\hat\TT_{r\setminus s})+\hat\TT_s(\hat\TT_s-\hat\TT_{s\setminus r})\right)\\
	&=\theta^2\hat\TT_{rs}\frac{\hat\TT_r+\hat\TT_s}{2}.\label{eqn:covAm}
\end{align}
Notice that the covariance formulations of Equations (\ref{eqn:covGm}) and (\ref{eqn:covAm}) are respectively based on the geometric mean and arithmetic mean of travel times $\hat\TT_r$ and $\hat\TT_s$. Note that although we base the covariance on the travel time and its overlap, it is also possible to define it based on overlap of distance or total systematic utility.  

\subsection{Desired choice model properties}\label{sec:rum:req}

In summary, the assumptions on the error terms and the utility formulation in Equation (\ref{eqn:U}) determine the expected value, variance and covariance of the utility. \begin{itemize}
\item The expected value of utility of route $r\in \rs$ is
\begin{align}
	\E (U_r)&=\E \left(V^0_r+\beta\TT_r+\rv_r\right),\\
	&=\E \left(V^0_r\right)+\E \left(\beta\TT_r\right)+\E \left(\rv_r\right),\\
	&=V^0_r+\beta\hat\TT_r+\E \left(\rv_r\right).
\end{align}
\item The variance of the utility of route $r\in \rs$  is 
\begin{align}
	\Var (U_r)&=\Var \left(V^0_r+\beta\TT_r+\rv_r\right),\\
	&=\Var \left(V^0_r\right)+\Var \left(\beta\TT_r\right)+\Var \left(\rv_r\right),\\
	&=\theta^2\hat\TT_r^2+\Var \left(\rv_r\right).
\end{align}
\item The covariance between routes $r,s\in \rs$,  $\Covar (U_r,U_s)=\Covar (\TT_r,\TT_s)$, equals either Equation (\ref{eqn:covGm}) or (\ref{eqn:covAm}).
\end{itemize}

Furthermore, when the route choice model is part of a transport network model, the following properties are desired:
\begin{itemize}
\item applicability to and compatibility between all types of (changing) networks, including those with OD pairs with short and long distances, and including those with changing differences and ratios between routes;
\item closed-form probabilities to avoid computationally expensive simulations; 
\item capture correlations due to overlap.
\end{itemize}

\section{Random Route Utility Maximization models with Generalized Multivariate  Extreme Value distributions}

In this section an overview is given of choice models that are applicable for routes based on random utility maximization (RUM) with closed-form expressions for the probabilities. This implies that MNP,  Mixed Logit, and models with an error correction component (e.g., Logit kernel) are not considered; we refer to \citep{Prashker2004,Frejinger2006_905} for an overview of these simulation based models. However, since the MNP model is capable of completely capturing both the analyst error and the random foreseen travel time (and thus heteroscedasticity and correlation)\footnote{Note that we use an adjusted MNP (compared to the literature) that is  in line with the (co)variances under the scaling postulate.}, we will use MNP simulated data as ground truth in the illustrative example in Section \ref{sec:network}. The existing models are all based on a single random variable, which should approximate both the analyst error and travel time distribution; Section \ref{sec:qualassess} answers to which degree this is possible.

In this section the family of distributions will be deduced from RUM models with separable systematic and random utility components (i.e., $U_r$ can be written as a function of $V_r$ and $\rv_r$). This leads to two types of Extreme Value (EV) models.  First consider the additive utility formula $U_r=V_r+\rv_r, \forall r\in\rs,$ where $\left(\rv_1,\ldots,\rv_{|\rs|}\right)$ is Multivariate Extreme Value (MEV) distributed \citep{McFadden1978}; these models fall in the EV type I (Gumbel) category. Second, consider the multiplicative utility formula $U_r=V_r\rv_r, \forall r\in\rs$ and assume that $\left(-\ln\rv_1,\ldots,-\ln\rv_{|\rs|}\right)$ is MEV distributed \citep{Galvez2002,Fosgerau2009494}; these models fall in the EV type III (reversed Weibull) category. The EV type II (Fr\'{e}chet) is not suitable for route choice since Fr\'{e}chet distributions have a lower bound on its support (i.e., there is a maximum on the route cost/travel time if such a model existed). The Generalized Multivariate Extreme Value (GMEV) models consist of both the additive and multiplicative MEV models. Note that in older literature the MEV models were also named GEV. GEV is a family of  univariate distributions consisting of Gumbel, Fr\'{e}chet and reversed Weibull. An MEV distribution is the joint distribution of multiple random variables with marginal distribution from one GEV type. MEV type I was discovered first by \citet{McFadden1978} under the name `The Generalized Extreme Value model'.

Other parametric approaches -- not based directly on separable RUM -- also exist. In \citep{Castillo2008373,Li2011,Nakayama2013753,Nakayama2015,Chikaraishi2015} utilities $U_r, \forall r\in\rs$ are independently Generalized Extreme Value (GEV) distributed. \citet{Li2011} captures even a larger class of models, with independent route utilities. \citet{Fosgerau20131} and \citet{Mattsson2014} provide the largest class of parametric approaches, where dependencies can be included as well. Note that all these models cannot be decomposed in systemtic and random components\footnote{Except for their specific cases where they fit in the \GMEV\ type.}. They directly parametrize the utility's probability distributions which is less explainable in behavioural terms. The exception is the derivation of the revised q-generalized logit model by \citet{Chikaraishi2015} from the q-product function. 

\subsection{The two types of Multivariate Extreme Value distributions}

MEV distributions are very flexible since they allow for correlation between their variables. Generating functions are the core of MEV distributions and unfortunately these are not easily interpretable, making the theory rather complex. It is well-known that MEV type I models capture the MNL, Paired Combinatorial Logit (PCL) and Link-Nested Logit (LNL) models, but it is less known that C-Logit, Path-Size Logit (PSL) and Path-Size Correction Logit belong to the MEV family. It is convenient that the latter models also fit in the generic framework since general results can be applied to them. All these mentioned `logit based' models have an additive utility formula, so first the MEV derivation for this additive case is presented. Thereafter, the multiplicative MEV models are derived, and we show how they both fit in the GMEV framework.

\subsubsection{Additive MEV models}
The additive utility formula is $U_r^\Aa=V_r+\rv^\Aa_r, \forall r\in\rs$. Define the probability of choosing route $r\in \rs$ as
\begin{align}
P_r^\Aa:= & \prob\left(U_r^\Aa\geq U_p^\Aa, \quad \forall p\neq r\right) \label{eqn:prob}\\
= & \prob\left(V_r+\rv_r^\Aa\geq V_p+\rv_p^\Aa, \quad \forall p\neq r\right).\label{eqn:Aprob}
\end{align}
Following \citet{McFadden1978}, assume that the joint distribution of $\left(\rv_1^\Aa,\ldots,\rv_{|\rs|}^\Aa\right)$ follows MEV type I distribution $\RV^\MEV$ with cumulative distribution function
\eqn{\label{eqn:addErrorDistr}
F_{\RV^\MEV}(x_1,\ldots,x_{|\rs|})=e^{-G(e^{-x_1},\dots,e^{-x_{|\rs|}})},
}
where \emph{generating function} $G:\mathbb{R}^{|\rs|}\rightarrow\mathbb{R}$ satisfies 
\begin{itemize}
\item $\scale$-homogeneity:
\eqn{G(a z_1,\ldots,a z_{|\rs|})=a^\scale G(z_1,\ldots,z_{|\rs|}), \quad \forall a>0\label{eqn:muhom}}
\item Limit property:
\eqn{\lim_{z_r\rightarrow\infty}G(z_1,\ldots,z_{|\rs|})=\infty, \quad\forall r\in\rs}
\item Alternating signs:
\eqn{
\frac{\partial^{|\hat\rs|}G(z_1,\ldots,z_{|\rs|})}{\prod_{r\in\hat{\rs}}\partial z_r}(-1)^{|\hat\rs|-1}\geq 0,\quad \forall \hat \rs \subseteq \rs, \hat \rs \neq \emptyset
.\label{eqn:altsign}}
\end{itemize}
The alternating signs property means that the sign of $G$ switches every time an additional distinct partial derivative is taken. The generating function needs to satisfy the properties in order for $F_{\RV^\MEV}$ to be a well-defined multivariate cumulative distribution function. With this utility formulation and distribution the following closed-form choice probabilities can be derived:
\begin{align}
P_r^\Aa =& \prob\left(U_r^\Aa\geq U_p^\Aa, \quad \forall p\neq r\right)\nonumber\\
=&\frac{y_r^\Aa G_r(y_1^\Aa,\ldots,y_{|\rs|}^\Aa)}{\scale G(y_1^\Aa,\ldots,y_{|\rs|}^\Aa)},\label{eqn:AMEVprob}
\end{align} 
where $G_r(x_1,\ldots,x_{|\rs|})=\frac{\partial G(x_1,\ldots,x_{|\rs|})}{\partial x_r}$ and $y_r^\Aa=e^{V_r},\forall r\in\rs$.
More details on the generating function and the derivation of the choice probabilities can be found in \citep{McFadden1978}. By using Euler's homogeneous function theorem Equation (\ref{eqn:AMEVprob}) can be rewritten as 
\begin{align}
	P_r^\Aa =&\frac{y_r^\Aa G_r(y_1^\Aa,\ldots,y_{|\rs|}^\Aa)}{\sum_{p\in\rs}y_p^\Aa G_p(y_1^\Aa,\ldots,y_{|\rs|}^\Aa)}\label{eqn:AMEVprob2}
	\\
	=&\frac{e^{V_r}G_r(e^{V_1},\ldots,e^{V_{|\rs|}})}{\sum_{p\in\rs}e^{V_p} G_p(e^{V_1},\ldots,e^{V_{|\rs|}})}\\
	=&\frac{e^{V_r+\ln G_r(e^{V_1},\ldots,e^{V_{|\rs|}})}}{\sum_{p\in\rs}e^{V_p+\ln G_p(e^{V_1},\ldots,e^{V_{|\rs|}})} }\label{eqn:AMEVprobeasy}
\end{align}
Since the denominator is independent of $r$ and the formulation is similar to that of MNL, the probabilities provided in Equation (\ref{eqn:AMEVprobeasy}) are the easiest to interpret. In the remainder of this paper the formulations of Equations (\ref{eqn:AMEVprob}) and (\ref{eqn:AMEVprob2}) will be used because in the multiplicative case only the definition of the $y_r^\Aa$'s will change, and these can then be substituted easily.

For a route set $\rs$ the combination of a generating function $G$ and a \emph{generating vector} $\y=(y_1,\ldots,y_{|\rs|})$ completely determines choice probabilities with Equation (\ref{eqn:AMEVprob}). In other words, the route set, the generating function and the generating vector completely specify a choice model. Therefore, we can write the choice probabilities as a function of the generating function and generating vector for all types of \GMEV\ models. For a choice model $X$ applied on route set $\rs$ with $\scale$-homogeneous generating function $G^X$ and generating vector $\y^X$, we can write the probability of choosing route $r$ as
\eqn{
P_r^X(G^X;\y^X):= \frac{y_r^XG_r^X(\y^X)}{\scale G^X(\y^X)},\quad \forall r\in \rs.\label{eqn:Gprobs}
}   

Note that for the additive utility formula the generating vector is always equal to the exponentials of systematic utilities. Therefore define $\yA:=(e^{V_1},\ldots,e^{V_{|\rs|}})$ as the \emph{additive generating vector}. This additive generating vector should be used for all MEV models with an additive utility form, thus the generating function then defines choice probabilities. 

\subsubsection{Multiplicative MEV models}

Before we present examples of the generating function, we first examine the multiplicative utility formula. \citet{Fosgerau2009494} are the first to explore this field. In their analysis, the multiplicative formulation is transformed into an equivalent additive formulation. The analysis presented here is slightly different because we allow more flexibility on the scale parameter. The models in \citep{Castillo2008373,Nakayama2013753,Kitthamkesorn2013} yield similar results for specific instances of generating functions, but the utility is not purely the product of systematic utility and an error term; instead, the utility is directly parametrized.

The multiplicative utility formula is 
\eqn{
U_r^\Mm=V_r\rv_r^\Mm,\quad \forall r\in\rs,
} 
where $V_r< 0$ and $\rv_r^\Mm \geq 0$. This domain of the systematic utility is more restrictive than for additive MEV; however, in the route choice context it is natural that route costs are valued negatively, yielding negative systematic utility. By applying a log-transformation to Equation (\ref{eqn:prob}) and some further algebra, we can derive the probability of choosing route $r\in\rs$ as:
%\begin{align}
%P_r^\Mm=&\prob\left(V_r\rv_r^\Mm\geq V_p\rv_p^\Mm, \quad \forall p\neq r\right)\label{eqn:MultAddStart}\\
%=&\prob\left(-\ln \left(-V_r\right)-\ln\left(\rv_r^\Mm\right)\geq-\ln\left( -V_p\right)-\ln\left(\rv_p^\Mm\right), \quad \forall p\neq r\right)\label{eqn:MultAdd}
%\end{align}
\begin{align}
	P_r^\Mm=&\prob\left(V_r\rv_r^\Mm\geq V_p\rv_p^\Mm, \quad \forall p\neq r\right)\label{eqn:MultAddStart}\\
	=&\prob\left(-V_r\rv_r^\Mm\leq -V_p\rv_p^\Mm, \quad \forall p\neq r\right)\\
	=&\prob\left(\ln \left(-V_r\rv_r^\Mm\right)\leq\ln\left( -V_p\rv_p^\Mm\right), \quad \forall p\neq r\right)\label{eqn:lntf}\\
	=&\prob\left(\ln \left(-V_r\right)+\ln\left(\rv_r^\Mm\right)\leq\ln\left( -V_p\right)+\ln\left(\rv_p^\Mm\right), \quad \forall p\neq r\right)\\
	=&\prob\left(-\ln \left(-V_r\right)-\ln\left(\rv_r^\Mm\right)\geq-\ln\left( -V_p\right)-\ln\left(\rv_p^\Mm\right), \quad \forall p\neq r\right)\label{eqn:MultAdd}
\end{align}
Equation (\ref{eqn:MultAdd}) has the same structure as the choice probabilities in Equation (\ref{eqn:Aprob}) of the additive choice model with systematic utilities $\tilde V_r = -\ln \left(-V_r\right), \forall r\in\rs$ and error terms $\tilde \rv_r^\Aa=-\ln\left(\rv_r^\Mm\right)$. Therefore, we can apply the theory of the additive MEV with these variable substitutions. So in the multiplicative case, $\left(-\ln\rv_1^\Mm,\ldots,-\ln\rv_{|\rs|}^\Mm\right)$ follows MEV type I distribution $\RV^\MEV$ (i.e., $(\rv_1^\Mm,\ldots,\rv_{|\rs|}^\Mm)$ is not similar to $\RV^\MEV$!), this distribution of $(\rv_1^\Mm,\ldots,\rv_{|\rs|}^\Mm)$ is called MEV type III. This is a multivariate distribution whose marginal distributions are of type reversed Weibull, and inherits the covariance structure from the generating function. So, it is the reversed Weibull equivalent of the model of \citet{McFadden1978}.  Note that all generating functions can be applied on both types. The \emph{multiplicative generating vector}, denoted with $\yM$, is derived from the additive generating vector and becomes 
\eqn{
\yM=\left(e^{\tilde V_1},\ldots,e^{\tilde V_{|\rs|}}\right)=\left(e^{-\ln \left(-V_1\right)},\ldots,e^{-\ln \left(-V_{|\rs|}\right)}\right)=\left(\frac{-1}{V_1},\ldots,\frac{-1}{V_{|\rs|}}\right).
}

The probabilities of models with multiplicative utility formulas can be derived from  Equation (\ref{eqn:Gprobs}) by simply using the multiplicative generating vector $\yM$. By applying Euler's homogeneous function theorem again we can get
\eqn{
P_r=\frac{\frac{-G_r(-1/V_1,\ldots,-1/V_{|\rs|})}{V_r}}{\sum_{p\in\rs}\frac{-G_p(-1/V_1,\ldots,-1/V_{|\rs|})}{V_p}}
,} but the form of Equation (\ref{eqn:Gprobs}) using $\yM$ is preferred. The main advantage of the multiplicative utility formula is that the utilities are automatically heteroscedastic. The standard deviation of the utility is proportional with the systematic utility and thus even the multiplicative equivalent of MNL is heteroscedastic. 

A well-known property of additive models is that the probabilities are invariant under addition of a constant to all systematic utilities. Therefore, only differences in utilities matter in the MNL model. Analogously, the multiplicative models are invariant under multiplying all systematic utilities with the same factor (mini-proof: If $\lambda V_r\rv_r^\Mm\geq \lambda V_p\rv_p^\Mm$ then $V_r\rv_r^\Mm\geq V_p\rv_p^\Mm$), and only ratios between utilities matter in the multiplicative equivalent of MNL. In Section \ref{sec:qualassess} its properties are further analysed.

\subsubsection{Hybrid approach}\label{sec:hybrid}

\citet{Xu2015} introduce a hybrid model directly derived from the logit and weibit choice probability functions. This model inherits characteristics of as well logit as weibit. The model uses two scale parameters, one from logit ($\scale^\Aa$) and one from weibit ($\scale^\Mm$), and has choice probabilities equal to
\begin{equation}\label{eqn:hybridXu}
P_r=\frac{e^{\scale^\Aa V_r}\left(-\frac{1}{V_r}\right)^{\scale^\Mm}}{\sum_{p\in\rs}e^{\scale^\Aa V_p}\left(-\frac{1}{V_p}\right)^{\scale^\Mm}},\quad \forall r\in \rs.
\end{equation}
Unfortunately, the model cannot be written as a RUM model with route $r$'s utility in the form of $U_r=V_r\times\rv^\Mm + \rv^\Aa$ with a additive error term $\rv^\Aa$ and multiplicative error term $\rv^\Mm$, which would have been a very neat solution. A straightforward hybrid approach in our framework would be to introduce generating vector $(-e^{V_1}/V_1,\ldots,-e^{V_{|\rs|}}/V_{|\rs|})$, i.e., the element-wise multiplication of $\yA$ and $\yM$. This will not lead to the approach of \citet{Xu2015}, since only one scale parameter will be introduced when the vector is applied on a generating function. 

However, the hybrid model of \citet{Xu2015} can be derived from either the \A\ or \M models from our framework by introducing an additional term in the systematic utility. In the additive form, the hybrid model of equation (\ref{eqn:hybridXu}) is obtained with 
\begin{equation}
U_r=V_r-\frac{\scale^\Mm}{\scale^\Aa}\ln(-x)+\rv^\Aa_r, \quad \forall r\in\rs,
\end{equation}
where a logarithmic term is inserted to simulate the multiplicative model \citep{Xu2015}. In our framework, this can be represented with generating vector 
\begin{equation}\label{eqn:yHybridAdd}
\y=\left(e^{V_1-\rho\ln(-V_1)},\ldots,e^{V_{|\rs|}-\rho\ln(-V_{|\rs|})}\right)=
\left(
	\frac
	{e^{V_1}}
	{\left(-V_1 \right)^\rho}
	,\ldots,
	\frac
	{e^{V_{|\rs|}}}
	{\left(-V_{|\rs|}\right)^\rho}
\right),
\end{equation}
where $\rho>0$ is a parameter representing the ratio between the scales. In the multiplicative form, the hybrid model of equation (\ref{eqn:hybridXu}) is obtained with 
\begin{equation}
U_r=V_r\times e^{-\frac{\scale^\Aa}{\scale^\Mm} V_r}\times \rv^\Mm_r, \quad \forall r\in\rs,
\end{equation}
where a exponential factor is inserted to simulate the additive model. In the framework again, this becomes generating vector 
\begin{equation}
\y=\left(\frac{-1}{V_1 e^{-\frac{V_1}{\rho}}},\ldots,\frac{-1}{V_{|\rs|} e^{-\frac{V_{|\rs|}}{\rho}}}\right)=
\left(
	\frac
	{e^{\frac{V_1}{\rho}}}
	{-V_1}
	,\ldots,
	\frac
	{e^{\frac{V_{|\rs|}}{\rho}}}
	{-V_{|\rs|}}
\right),
\end{equation}
with the same interpretation for $\rho$ as above. The difference between the two generating vectors lies at the estimated scale parameter. With the latter formulation $\scale^\Mm$ is estimated directly from the model and $\scale^\Aa$ is derived from $\rho$. While with equation (\ref{eqn:yHybridAdd}) $\scale^\Aa$ is estimated directly and $\scale^\Mm$ is derived from $\rho$. There is no difference in model outcomes and/or flexibility. Using any of the two approaches, it is possible to derive hybrid path-size, paired combinatorial, and link-nested models (see next section). The covariance structure provided in the latter two model instances are then only applied on the main error structure (i.e., dependent on which of the two generating functions is chosen).
%\escom{Should we move to section to the end of the paper? At the discussion? Or somewhere else?} 

\subsection{Generating functions and model instances}\label{sec:instances}
To create instances of the family of GMEV route choice models, the generating function has to be specified. Together with one of the two derived generating vectors, Equation (\ref{eqn:Gprobs}) will provide the closed-form choice probabilities of the choice model. Any function that satisfies Equations (\ref{eqn:muhom})-(\ref{eqn:altsign}) is a generating function. The cross-nested logit uses one of the most general forms of the generating function and is analysed by \citet{Bierlaire2006,Papola2004}. All presented generating functions here are also special cases of cross-nested logit models.

To address overlap on links, the presented model instances assume that each link is associated with one cost, namely systematic utility. This assumption is made for clarity reasons, differentiation into link length, travel and other cost can easily be done, also for the various factors and coefficients introduced below. Let $\LL$ be a set of links and let the systematic utility of link $l\in \LL$ be $V_l$. For each route $r\in \rs$, $\LL_r$ is the set of links of which $r$ consists and the systematic route utility equals $V_r=\sum_{l\in\LL_r}V_l$. 
 
\subsubsection{Multinomial}
The well-known MNL model \citep[see ][]{McFadden74} is based on multinomial generating function
\eqn{
\GMNL(\z):=\sum_{r\in\rs}z_r^\scale,
}
where $\scale>0$ is the scale parameter. The traditional (additive) MNL model, denoted \A\MNL, has the following simple route choice probabilities,
\eqn{ \label{eqn:AMNLprobs}
P_r^{\A\MNL}(\GMNL;\yA)= \frac{y^\text{A}_r\GMNL_r(\yA)}{\scale \GMNL(\yA)}=\frac{e^{\scale V_r}}{\sum_{p\in\rs}e^{\scale V_p}},\quad \forall r\in \rs.
}  
The multiplicative counterpart is the multinomial weibit model (or Kirchhoff distribution), denoted with \M\MNL\ and has route choice probabilities
\eqn{
P_r^{\M\MNL}(\GMNL;\yM)= \frac{y^\text{M}_r\GMNL_r(\yM)}{\scale \GMNL(\yM)}=\frac{\left(-\frac{1}{V_r}\right)^\scale}{\sum_{p\in\rs}\left(-\frac{1}{V_p}\right)^\scale},\quad \forall r\in \rs.
}
\subsubsection{Path-Size}
We will first introduce a Path-Size (PS) generating function to define a specific MEV type, and then we will show that models proposed in the literature are identical. For this PS generating function the path-size factor $\ps_r$ is required for each route $r\in\rs$; this factor depicts the amount of overlap with other routes. Define PS generating function
\eqn{
\GPSL(\z):=\sum_{r\in\rs}\ps_r^\pspar z_r^\scale,
}
where $\scale>0$ is the scale parameter and $\pspar$ is the path-size parameter. The additive Path-Size model, denoted \A\PSL, has choice probabilities 
\eqn{
P_r^{\A\PSL}(\GPSL;\yA)= \frac{y^\text{A}_r\GPSL_r(\yA)}{\scale \GPSL(\yA)}=\frac{\ps_r^\pspar e^{\scale V_r}}{\sum_{p\in\rs}\ps_p^\pspar e^{\scale V_p}},\quad \forall r\in \rs.\label{eqn:PA-PSL}
}  

The path-size-like models in the literature have a logarithmic term included in the utility. The MNL model with scale $\scale$ and systematic route utility $V_r+\gamma\ln x_r$ leads to the same choice probabilities as the additive MEV model with generating function $\GPSL(\z)$, generating vector $\y=(e^{V_1},\dots,e^{V_|\rs|})$, path-size factors $\ps_r=x_r, \forall r\in\rs$, and path-size parameter $\pspar=\scale\gamma$. So, C-Logit, Path-Size Logit, and Path-Size Correction Logit can be translated into a \A\PSL\ model. C-Logit presented by \citet{Cascetta1996} adds the term $-\beta^\text{CF}\text{CF}_r$ to the utility, with parameter $\beta^\text{CF}$ and commonality factor $\text{CF}_r$ for each route $r$. This is equivalent to \A\PSL\ with $\ps_r=\text{CF}_r$ and $\pspar=-\scale\beta^\text{CF}$. Path-Size Logit presented by \citet{ben-akiva1999} adds the term $\beta^\text{PS} {\tilde \ps}_r$ to the utility, with parameter $\beta^\text{PS}$ and path-size factor ${\tilde \ps}_r$ for each route $r$. This is equivalent to \A\PSL\ with $\ps_r={\tilde \ps}_r$ and $\pspar=\scale\beta^\text{PS}$. Path-Size Correction Logit \citep[see ][]{bovy2008} is equal in the same way. For those familiar with general cross-nested logit: note that \GPSL\ is an instance of cross-nested logit with parametrized non-normalized inclusion factors with a single nest per alternative.

Different choices for $\ps_r$ are compared in \citep{Frejinger2006_905}; they conclude that the best formulation for the path-size factor is
\begin{equation}\label{eqn:pathsize}
\ps_r=\frac{\sum_{l\in \LL_r}\frac{V_l}{\#_l}}{\sum_{l\in \LL_r}V_l}=\frac{\sum_{l\in \LL_r}\frac{V_l}{\#_l}}{V_r},
\end{equation}
where $\#_l=\left|\{r\in \rs| l\in \LL_r\}\right|$ is the number of routes using link $l$. If $\ps_r=1$ then $r$ has no overlap with any other route. If $\ps_r\rightarrow 0$ then $r$ shares each link with many other routes. 

The multiplicative counterpart of \A\PSL\ is denoted with \M\PSL\ and has route choice probabilities
\eqn{
P_r^{\M\PSL}(\GPSL;\yM)= \frac{y^\text{M}_r\GPSL_r(\yM)}{\scale \GPSL(\yM)}=\frac{\ps_r^\pspar\left(-\frac{1}{V_r}\right)^\scale}{\sum_{p\in\rs}\ps_p^\pspar\left(-\frac{1}{V_p}\right)^\scale},\quad \forall r\in \rs.
}
Path-size is also added to MNW by \citet{Kitthamkesorn2013,Kitthamkesorn2014}, but they assume that $\pspar=1$, which thus leads to a specific instance of \M\PSL.

\subsubsection{Paired Combinatorial}
The (additive) Paired Combinatorial Logit (PCL) model as presented and analysed by \citet{Chu1989,Koppelman200075,pravinvongvuth2005,Prashker1998,gliebe1999route,Chen2003} captures correlation between each pair of routes. For each route couple a nest is defined with a fixed nest specific scale. This scale is determined by the similarity index, denoted with $\cf_{rp}$, for all $r\neq p\in\rs$. Define paired combinatorial generating function
\eqn{
\GPCL(\z):=\sum_{r\in\rs}\sum_{p\in\rs\setminus\{r\}}
\left(
z_r^\frac{\scale}{1-\cf_{rp}}+z_p^\frac{\scale}{1-\cf_{rp}}
\right)^{1-\cf_{rp}},
}
where $\scale>0$ is the scale parameter. Multiple definitions for the similarity index exist, but the most common one is
\begin{equation}
\cf_{rp}:=\frac{\sum_{l\in \LL_r\cap\LL_p}V_l}{\sqrt{V_rV_p}}, \quad \forall r,p\in\rs.
\end{equation}

The traditional additive PCL model, denoted \A\PCL, has route choice probabilities 
\eqn{
P_r^{\A\PCL}(\GPCL;\yA)= \frac{y^\text{A}_r\GPCL_r(\yA)}{\scale \GPCL(\yA)}=
\frac{
 \sum_{p\in\rs\setminus\{r\}} e^\frac{\scale V_r}{1 - \cf_{rp}} \left(e^\frac{\scale V_r}{1 - \cf_{rp}}+e^\frac{\scale V_p}{1 - \cf_{rp}}\right)^{-\cf_{rp}}
}{
 \sum_{r'\in\rs}\sum_{p\in\rs\setminus\{r'\}}\left(e^\frac{\scale V_{r'}}{1 - \cf_{r'p}}+e^\frac{\scale V_p}{1 - \cf_{r'p}}\right)^{1-\cf_{r'p}}
}
,\quad \forall r\in \rs.
}  
The multiplicative counterpart is denoted with \M\PCL\ and has route choice probabilities
\eqn{
P_r^{\M\PCL}(\GPCL;\yM)= \frac{y^\text{M}_r\GPCL_r(\yM)}{\scale \GPCL(\yM)}=
\frac{
 \sum_{p\in\rs\setminus\{r\}} \left(-\frac{1}{V_r}\right)^\frac{\scale }{1 - \cf_{rp}} \left(\left(-\frac{1}{V_r}\right)^\frac{\scale}{1 - \cf_{rp}}+\left(-\frac{1}{V_p}\right)^\frac{\scale}{1 - \cf_{rp}}\right)^{-\cf_{rp}}
}{
 \sum_{r'\in\rs}\sum_{p\in\rs\setminus\{r'\}}\left(\left(-\frac{1}{V_{r'}}\right)^\frac{\scale}{1 - \cf_{r'p}}+\left(-\frac{1}{V_p}\right)^\frac{\scale}{1 - \cf_{r'p}}\right)^{1-\cf_{r'p}}
}
,\quad \forall r\in \rs.
}
\subsubsection{Link-Nested}
\citet{vovsha1998} are the first who applied the Cross-Nested Logit model to route choice such that links represent nests. In this (additive) Link-Nested Logit (LNL) a nest is created for each link and all routes that use the link are included in the nest. The inclusion coefficient of each route in each link (or nest) is denoted with $\ic_{lr}$ for all links $l\in\LL$ and routes $r\in\rs$. Define link-nested generating function
\eqn{
\GLNL(\z):=\sum_{l\in\LL} \left(\sum_{r\in\rs}\ic_{lr}z_r^{\scale_l}\right) ^\frac{\scale}{\scale_l},
}
where $\scale>0$ is the scale parameter and $\scale_l>0, \forall l\in\LL$ are the link specific scale parameters. In general, it is difficult to estimate all link specific scale parameters in large networks. More on this topic can be found in \citep{Bierlaire2006}, where also the necessity of normalizing the inclusion coefficients is discussed. For this study the following normalized inclusion coefficients are applied
\begin{equation}
\ic_{lr}=\begin{cases}
\frac{V_l}{V_r} & \text{if } l\in\LL_r\\
0 & \text{otherwise}
\end{cases}.
\end{equation}
 
 The traditional additive LNL model, denoted \A\LNL, has  route choice probabilities 
\eqn{
P_r^{\A\LNL}(\GLNL;\yA)= \frac{y^\text{A}_r\GLNL_r(\yA)}{\scale \GLNL(\yA)}=
\frac{
 \sum_{l\in\LL} \ic_{lr}e^{\scale_l V_r} \left(\sum_{p\in\rs}\ic_{lp}e^{\scale_lV_p}\right)^{\frac{\scale}{\scale_l}-1}
}{
 \sum_{l\in\LL}\left(\sum_{p\in\rs}\ic_{lp}e^{\scale_lV_p}\right)^\frac{\scale}{\scale_l}
}
,\quad \forall r\in \rs.
}  
The multiplicative counterpart is \M\LNL\ with route choice probabilities
\eqn{
P_r^{\M\LNL}(\GLNL;\yM)= \frac{y^\text{M}_r\GLNL_r(\yM)}{\scale \GLNL(\yM)}=
\frac{
 \sum_{l\in\LL} \ic_{lr}\left(-\frac{1}{V_r}\right)^{\scale_l} \left(\sum_{p\in\rs}\ic_{lp}\left(-\frac{1}{V_p}\right)^{\scale_l}\right)^{\frac{\scale}{\scale_l}-1}
}{
 \sum_{l\in\LL}\left(\sum_{p\in\rs}\ic_{lp}\left(-\frac{1}{V_p}\right)^{\scale_l}\right)^\frac{\scale}{\scale_l}
}
,\quad \forall r\in \rs.
}

\subsubsection{Joint Network}\label{sec:jointnetwork}

Recently, \citet{Papola2013CACIE} introduced a new generating function based on an generalization of cross-nested generating functions to networks \citep[see ][]{Daly2006,Newman2008}. Because their approach is based on links instead of routes, the generating function defines dependencies between links instead of routes. This is fundamentally different from the GMEV models presented in this work that is based on a route set and based on the error structure of random utility per route. The joint network formulation of \citet{Papola2013CACIE} is not a member of the \GMEV\ models for route choice in this study since it cannot handle generic route sets as input. Their method has an implicit enumeration of routes. This has the advantage that the problem of choice set generation is dealt with internally. Unfortunately, this does not allow a utility formulation per route. Furthermore, the link based approach also delimits the routes that will be found. For their model specifically, scale parameters of `nodes' depend on the shortest path to the destination. Their definition requires that a route advances through nodes for which the shortest path towards the destination decreases in every step. In real networks, the remaining shortest path will often increase when you deviate from the shortest path.  This makes the model restrictive. Although a natural connection between a road network and a (recursive) network based generating function seems to be promising, \citet{Marzano2014} suggests that the achievable covariance structure is not more advanced than that of cross-nested generating functions.

\section{Qualitative Assessment of the Models}\label{sec:qualassess}

This section analyses the utility distribution of the \GMEV\ models and assesses them qualitatively based on the desired properties presented in Section \ref{sec:rum:req}. The desired route utilities have two random variables (i.e., the random foreseen travel time and the analyst error) per route, but the models contain only one (\GMEV -distributed) variable per route. We show that this variable can only resemble the analyst error in the additive models, and that it can resemble both errors in the multiplicative models, however, with completely correlated foreseen travel times and analyst errors.  

\subsection{Utility distribution}

The results in this section for the additive models summarize the findings of \citet{McFadden1978,Daly2006}, while the derivations for the multiplicative models are slightly different from those in \citep{Fosgerau2009494} that uses two scales. From the additive utility formulation and Equation (\ref{eqn:addErrorDistr}) follows the multivariate utility distribution of the additive models:
%\begin{align}
%F_{U^\Aa}(x_1,\ldots,x_{|\rs|})=&\Pr\left(V_r +\rv_r^\Aa\leq x_r,\quad\forall r\in \rs\right)\\
%=&F_{\RV^\MEV}\left(x_1-V_1,\ldots,x_{|\rs|}-V_{|\rs|}\right)\\
%=&e^{-G(e^{-x_1+V_1},\dots,e^{-x_{|\rs|}+V_{|\rs|}})}
%\end{align}
\begin{align}
	F_{U^\Aa}(x_1,\ldots,x_{|\rs|})=&\Pr\left(U_r^\Aa\leq x_r,\quad\forall r\in \rs\right)\\
	=&\Pr\left(V_r +\rv_r^\Aa\leq x_r,\quad\forall r\in \rs\right)\\
	=&\Pr\left(\rv_r^\Aa\leq x_r-V_r ,\quad\forall r\in \rs\right)\\
	=&F_{\RV^\MEV}\left(x_1-V_1,\ldots,x_{|\rs|}-V_{|\rs|}\right)\\
	=&e^{-G(e^{-x_1+V_1},\dots,e^{-x_{|\rs|}+V_{|\rs|}})}
\end{align}
With some algebra, the multivariate utility distribution of the multiplicative models is written as
%\begin{align}
%F_{U^\Mm}(x_1,\ldots,x_{|\rs|})=&\Pr\left(V_r \rv_r^\Mm\leq x_r,\quad\forall r\in \rs\right)\\
%=&\Pr\left(-\ln \left(\rv_r^\Mm\right)\leq\ln \left(-V_r \right)- \ln\left(-x_r\right),\quad\forall r\in \rs\right)\\
%=&F_{\RV^\MEV}\left(\ln\left(-V_1 \right)- \ln\left(-x_1\right),\ldots,\ln\left(-V_{|\rs|} \right)- \ln\left(-x_{|\rs|}\right)\right)\\
%=&e^{-G\left(\frac{x_1}{V_1},\dots,\frac{x_{|\rs|}}{V_{|\rs|}}\right)}
%\end{align}
\begin{align}
	F_{U^\Mm}(x_1,\ldots,x_{|\rs|})=&\Pr\left(U_r^\Mm\leq x_r,\quad\forall r\in \rs\right)\\
	=&\Pr\left(V_r \rv_r^\Mm\leq x_r,\quad\forall r\in \rs\right)\\
	=&\Pr\left(-V_r \rv_r^\Mm\geq -x_r,\quad\forall r\in \rs\right)\\
	=&\Pr\left(\ln \left(-V_r \rv_r^\Mm\right)\geq \ln\left(-x_r\right),\quad\forall r\in \rs\right)\\
	=&\Pr\left(\ln \left(-V_r \right)+\ln \left(\rv_r^\Mm\right)\geq \ln\left(-x_r\right),\quad\forall r\in \rs\right)\\
	=&\Pr\left(-\ln \left(\rv_r^\Mm\right)\leq\ln \left(-V_r \right)- \ln\left(-x_r\right),\quad\forall r\in \rs\right)\\
	=&F_{\RV^\MEV}\left(\ln\left(-V_1 \right)- \ln\left(-x_1\right),\ldots,\ln\left(-V_{|\rs|} \right)- \ln\left(-x_{|\rs|}\right)\right)\\
	=&F_{\RV^\MEV}\left(\ln\frac{V_1}{x_1},\ldots,\ln\frac{V_{|\rs|}}{x_{|\rs|}}\right)\\
	=&e^{-G\left(e^{-\ln\frac{V_1}{x_1}},\dots,e^{-\ln\frac{V_{|\rs|}}{x_{|\rs|}}}\right)}\\
	=&e^{-G\left(\frac{x_1}{V_1},\dots,\frac{x_{|\rs|}}{V_{|\rs|}}\right)}
\end{align}
Then the marginal distribution for the utility of route $r\in\rs$ is 
%\begin{align}
%F_{U_r^\Aa}(x_r)=&\lim_{\left\{x_p\rightarrow \infty\right\}_{p\neq r}}F_{U^\Aa}\left(x_1,\ldots,x_{|\rs|}\right)\\
%=&e^{-G\left(e^{-x_r+V_r}\one_r\right)}\\
%=&e^{-e^{\scale \left(-x_r+V_r+\frac{\ln G\left(\one_r\right)}{\scale}\right)}},
%\end{align}
\begin{align}
	F_{U_r^\Aa}(x_r)=&\lim_{\left\{x_p\rightarrow \infty\right\}_{p\neq r}}F_{U^\Aa}\left(x_1,\ldots,x_{|\rs|}\right)\\
	=&e^{-G\left(e^{-x_r+V_r}\one_r\right)}\\
	=&e^{-e^{\scale (-x_r+V_r)}G\left(\one_r\right)}\\
	=&e^{-e^{\scale (-x_r+V_r)+\ln G\left(\one_r\right)}}\\
	=&e^{-e^{\scale \left(-x_r+V_r+\frac{\ln G\left(\one_r\right)}{\scale}\right)}},
\end{align}
for the additive case, which can be identified as a Gumbel distribution, and with notation $\one_r:=\left(0,\ldots,0,1,0,\ldots,0\right)$ (1 on the $r$-th position).\footnote{Note that for $\scale=1$ these are reversed exponential distributions.} For the multiplicative models the marginal distribution function of route $r$'s utility is 
%\begin{align}
%F_{U_r^\Mm}(x_r)=&\lim_{\left\{x_p\rightarrow 0\right\}_{p\neq r}}F_{U^\Mm}\left(x_1,\ldots,x_{|\rs|}\right)\\
%=&e^{-G\left(\frac{x_r}{V_r}\one_r\right)}\\
%=&e^{-\left(\frac{x_r}{V_rG\left(\one_r\right)^{\frac{-1}{\scale}}}\right)^\scale},
%\end{align}
\begin{align}
	F_{U_r^\Mm}(x_r)=&\lim_{\left\{x_p\rightarrow 0\right\}_{p\neq r}}F_{U^\Mm}\left(x_1,\ldots,x_{|\rs|}\right)\\
	=&e^{-G\left(\frac{x_r}{V_r}\one_r\right)}\\
	=&e^{-\left(\frac{x_r}{V_r}\right)^\scale G\left(\one_r\right)}\\
	=&e^{-\left(\frac{x_r}{V_rG\left(\one_r\right)^{\frac{-1}{\scale}}}\right)^\scale},
\end{align}
which can be identified as a reversed Weibull distribution, where reversed means that $-U_r^\Mm$ is Weibull distributed .

The expected maximum utility is an important measure of choice models; it can be used to formulate the corresponding user equilibrium. Furthermore, it is input to define duality gaps. Denote $U^{*\Aa}$ as the maximum utility for the additive models, and $U^{*\Mm}$ as the maximum utility for the multiplicative models. Their distributions can be written as
%\begin{align}
%F_{U^{*\Aa}}(x)=&F_{U^\Aa}(x,\ldots,x)\\
%=&e^{-G\left(e^{-x+V_1},\dots,e^{-x+V_{|\rs|}}\right)}\\
%=&e^{-e^{\scale \left(-x+\ln G\left(e^{V_1},\dots,e^{V_{|\rs|}}\right)/\scale\right)}}, \text{ and}\\
%F_{U^{*\Mm}}(x)=&F_{U^\Mm}(x,\ldots,x)\\
%=&e^{-G\left(\frac{x}{V_1},\dots,\frac{x}{V_{|\rs|}}\right)}\\
%=&e^{-\left(\frac{x}{-G\left(\frac{-1}{V_1},\dots,\frac{-1}{V_{|\rs|}}\right)^{\frac{-1}{\scale}}}\right)^\scale},
%\end{align} 
\begin{align}
	F_{U^{*\Aa}}(x)=&F_{U^\Aa}(x,\ldots,x)\\
	=&e^{-G\left(e^{-x+V_1},\dots,e^{-x+V_{|\rs|}}\right)}\\
	=&e^{-e^{-\scale x}G\left(e^{V_1},\dots,e^{V_{|\rs|}}\right)}\\
	=&e^{-e^{\scale \left(-x+\ln G\left(e^{V_1},\dots,e^{V_{|\rs|}}\right)/\scale\right)}}, \text{ and}\\
	F_{U^{*\Mm}}(x)=&F_{U^\Mm}(x,\ldots,x)\\
	=&e^{-G\left(\frac{x}{V_1},\dots,\frac{x}{V_{|\rs|}}\right)}\\
	=&e^{-(-x)^\scale G\left(\frac{-1}{V_1},\dots,\frac{-1}{V_{|\rs|}}\right)}\\
	=&e^{-\left(\frac{x}{-G\left(\frac{-1}{V_1},\dots,\frac{-1}{V_{|\rs|}}\right)^{\frac{-1}{\scale}}}\right)^\scale},
\end{align} 
which can be identified as Gumbel and reversed Weibull distributions again. 

%\begin{align}
%F_{U_r^\Aa,U_s^\Aa}(x_r,x_s)=&\lim_{\left\{x_p\rightarrow \infty\right\}_{\substack{p\neq r\\ p\neq s}}}F_{U^\Aa}\left(x_1,\ldots,x_{|\rs|}\right)\\
%=&e^{-G(0,\ldots,0,e^{-x_r+V_r},0,\ldots,0,e^{-x_s+V_s},0,\ldots,0)}
%\end{align}
%
%\begin{align}
%F_{U_r^\Mm,U_s^\Mm}(x_r,x_s)=&\lim_{\left\{x_p\rightarrow 0\right\}_{\substack{p\neq r\\ p\neq s}}}F_{U^\Mm}\left(x_1,\ldots,x_{|\rs|}\right)\\
%=&e^{-G\left(0,\ldots,0,\frac{x_r}{V_r},0,\ldots,0,\frac{x_s}{V_s},0,\ldots,0\right)}
%\end{align}
%\begin{align}
%\Var(U_r^\Aa)=\int_{\xi=-\infty}^\infty \xi^2 f_{U_r^\Aa}(\xi)\diff \xi-\left(\int_{\xi=-\infty}^\infty \xi f_{U_r^\Aa}(\xi)\diff \xi\right)^2, \quad \forall r\in \rs,
%\end{align}
%where $f_{U_r^\Aa}(x)=\diff F_{U_r^\Aa}(x)/\diff x$ is the probability density function of $U_r^\Aa$
%
%\begin{align}
%\Covar(U_r^\Aa,U_s^\Aa)=\int_{\xi_r=-\infty}^\infty\int_{\xi_s=-\infty}^\infty F_{U_r^\Aa,U_s^\Aa}(\xi_r,\xi_s)-F_{U_r^\Aa}(\xi_r)F_{U_s^\Aa}(\xi_s)\diff\xi_s\diff \xi_r, \quad \forall r,s\in \rs,
%\end{align}
%
%\begin{align}
%\Var(U_r^\Mm)=\int_{\xi=-\infty}^0 \xi^2 f_{U_r^\Mm}(\xi)\diff \xi-\left(\int_{\xi=-\infty}^0 \xi f_{U_r^\Mm}(\xi)\diff \xi\right)^2, \quad \forall r\in \rs, 
%\end{align}
%where $f_{U_r^\Mm}(x)=\diff F_{U_r^\Mm}(x)/\diff x$ is the probability density function of $U_r^\Mm$.
%
%\begin{align}
%\Covar(U_r^\Mm,U_s^\Mm)=\int_{\xi_r=-\infty}^0\int_{\xi_s=-\infty}^0 F_{U_r^\Mm,U_s^\Mm}(\xi_r,\xi_s)-F_{U_r^\Mm}(\xi_r)F_{U_s^\Mm}(\xi_s)\diff\xi_s\diff \xi_r, \quad \forall r,s\in \rs,
%\end{align}

\begin{table}
\caption{Expected value and variance of the desired, multiplicative and additive utility formulations}
	\centering
\begin{tabular}{l|lll}
\toprule  & Desired & Additive MEV & Multiplicative MEV \\ 
  & & {(\small $\gamma$ is Euler's constant)} & {(\small $\Gamma(\cdot)$ is the Gamma function)} \\ 
\midrule Utility $U_r$ & $V^0_r+\beta\TT_r+\rv_r$ & $V_r+\rv_r^\Aa=V^0_r+\beta\hat\TT_r+\rv_r^\Aa$ & $V_r\times\rv_r^\Mm=\left(V^0_r+\beta\hat\TT_r\right)\times\rv_r^\Mm$ \\ 
Expected value  & $\underbrace{V^0_r+\beta\hat\TT_r+\E \left(\rv_r\right)}_{\text{affine in }\hat\TT_r}$ & $\displaystyle \underbrace{V^0_r+\beta\hat\TT_r+\frac{\ln G(\one_r)+\gamma}{\scale}}_{\text{affine in }\hat\TT_r}$ & $\displaystyle \underbrace{\frac{V^0_r+\beta\hat\TT_r}{G(\one_r)^{1/\scale}}\Gamma\left(1+\frac{1}{\scale}\right)}_{\text{affine in }\hat\TT_r}$\\ 
Variance  & $\underbrace{\theta^2\hat\TT_r^2+\Var \left(\rv_r\right)}_{\text{quadratic in }\hat\TT_r}$ & $\displaystyle\underbrace{\frac{\pi^2}{6}\left(\frac{1}{\scale}\right)^2}_{\text{constant}}$ & $\displaystyle\underbrace{\frac{\left(V^0_r+\beta\hat\TT_r\right)^2}{G(\one_r)^{2/\scale}}\left(\Gamma\left(1+\frac{2}{\scale}\right)-\Gamma\left(1+\frac{1}{\scale}\right)^2\right)}_{{\text{quadratic in }V^0_r+\beta\hat\TT_r}}$ \\  
Expected $\max_{r\in\rs} U_r$  & Not closed form & $\displaystyle \frac{\ln G\left(e^{V_1},\dots,e^{V_{|\rs|}}\right)+\gamma}{\scale}$ & $\displaystyle -G\left(\frac{-1}{V_1},\dots,\frac{-1}{V_{|\rs|}}\right)^{\frac{-1}{\scale}}\Gamma\left(1+\frac{1}{\scale}\right)$ \\
\bottomrule 
\end{tabular} 
\label{tab:utility}
\end{table}

Table \ref{tab:utility} shows the desired expected value, variance, and expected maximum utility together with the actual expected value, variance and expected maximum utility of the models. It is based on the identification of the marginal distributions as Gumbel and reversed Weibull, and on systematic utility specification $V_r=V^0_r+\beta\TT_r$. Neither the additive nor the multiplicative formulation completely coincide with the desired result. 

The expected value is not a problem for any model; any value can be achieved after normalization and identification, and the location of utility is not decisive for choice probabilities. On the other hand, the variance of the additive models is constant and the variance of the multiplicative models is affected by $V_r^0$, and for the latter $\Var(\rv_r)$ is not directly represented. However, as the next section describes, the constant in the systematic utility of the multiplicative model does not have to be normalized. This constant will return as a constant term in the variance, and can thus resemble $\Var(\rv_r)$. Despite that the multiplicative models' variances can capture the variances from the foreseen travel time distribution and analyst error simultaneously, they are both represented by $\rv_r$ and thus completely dependent. This is not in line with the desired independence.    

Closed form formulations of the covariance of GMEV models are not known. \citet{Marzano20131,Marzano2014} present a more tractable expression and a method to calculate them for additive MEV. The covariances of the additive models are based on the constant variance, and will therefore not change when the systemic utility changes. However, the systematic utility enters the variance quadratically in the multiplicative model, and therefore its covariance will also in- and decrease together with it. Thus, the multiplicative models capture the covariances better.  

The expected maximum utility for the \A\MNL\ model indeed return in the well-known log-sum that can by used for an hierarchical derivation of the other models. This was first identified by \citet{ben1973structure}.

\subsection{Stochastic user equilibrium formulation}

The derived \GMEV\ route choice models all have an equivalent stochastic user equilibrium formulation. Mathematical programming formulation are known for all logit based models \citep{Fisk1980,Bekhor1999}, for the \M\MNL\ and (simplified) \M\PSL\ models \citep{Kitthamkesorn2013,Kitthamkesorn2014}, for the q-generalized logit model \citep{Chikaraishi2015}, and for the hybrid logit-weibit model \citep{Xu2015}. It is straightforward to derive the variational inequality formulation -- that we discuss -- from a mathematical programming formulation. \citet{chen1999dynamic,Guo2010} describe a stochastic user equilibrium for the MNL model using a variational inequality. \citet{Zhou2012} provide two variational inequality formulations for the C-logit model. Their first formulation is specific for C-logit and a special case of the formulation we present below. Their second formulation is more generic and works with any choice probability formula, but it returns different generalised costs than the approach we present below.

Since the choice probabilities of each \GMEV\ model can be written solely in terms of the generating function and generating vector, this also holds for the corresponding variational inequality  formulation. Denote the demand with $\demand$, the flow for route $r\in\rs$ as $\flow_r$, and $\flowvec$ as the vector of flows. Since, amongst other attributes, the travel time in the systematic utility depends on the flow, write that the generating vector now depends on $\y(\flowvec)$\footnote{When the generating vector is substituted, just replace $V_r$ with $V_r(\flowvec)$ to denote dependency of attributes on flow.}. For clarity, but without loss of generality, only one route set (i.e., one OD pair) is considered. 

Consider the following variational inequality (VI) formulation; find equilibrium flow $\flowvec^*=(\flow_1^*,\ldots,\flow_{|\rs|}^*)$, such that
\begin{equation}
\begin{split}
\sum_{r\in\rs} \left(-\ln\left(y_r(\flowvec^*)G_r(\y(\flowvec^*))\right)+\ln(\flow^*_r)\right)(\flow_r-\flow_r^*)\geq 0,\quad \forall \flowvec\in\flowfeas,\\
\text{where }\flowfeas=\left\{\flowvec\in\R^{|\rs|}\middle|f_r>0,\sum_{r\in\rs}\flow_r=\demand\right\}.
\end{split}\label{eqn:suevi}
\end{equation} 
Define the generalized stochastic cost for route $r\in\rs$ under flow $\flowvec$ as $c_r(\flowvec):=-\ln\left(y_r(\flowvec)G_r(\y(\flowvec))\right)+\ln(\flow_r)$\footnote{Note that this cost is unit-less. For MNL, it is possible to normalize it such that the systemic utility appears as a term.}. The corresponding Karush-Kuhn-Tucker system \citep[see e.g.,][Proposition 1.2.1]{facchinei2003finite} is then to find multipliers $\kkta$ and $\kktb_1,\ldots,\kktb_{|\rs|}$ for which
\begin{align}
0&=-\ln\left(y_r(\flowvec^*)G_r(\y(\flowvec^*))\right)+\ln(\flow_r^*)+\kkta-\kktb_r,\quad \forall r\in\rs,\label{eqn:kkt1}\\
\sum_{r\in\rs}\flow_r&=\demand,\\
0&\leq\kktb_r, \quad \forall r\in\rs, \label{eqn:kkt2}\\
0&<\flow_r, \quad \forall r\in\rs, \text{ and}\label{eqn:kkt3}\\
0&=\sum_{r\in\rs}\kktb_r\flow_r\label{eqn:kkt4}
\end{align} 
hold. Equations (\ref{eqn:kkt2}), (\ref{eqn:kkt3}) and (\ref{eqn:kkt4}) imply that $\kktb_r=0$ for all $r\in\R$. Rewriting equation (\ref{eqn:kkt1}) then gives 
\begin{align}
\ln\left(y_r(\flowvec^*)G_r(\y(\flowvec^*))\right)&=\ln(\flow_r^*)+\kkta,\quad \forall r\in\rs\\
y_r(\flowvec^*)G_r(\y(\flowvec^*))&=\flow^*_re^\kkta,\quad \forall r\in\rs.\label{eqn:sueroute}
\end{align}
Summing over routes gives
\begin{equation}
\sum_{r\in\rs}y_r(\flowvec^*)G_r(\y(\flowvec^*))=\demand e^\kkta,\quad \forall r\in\rs.\label{eqn:suesum}
\end{equation}
Finally, dividing equation (\ref{eqn:sueroute}) by equation (\ref{eqn:suesum}) gives
\begin{equation}
\frac{y_r(\flowvec^*)G_r(\y(\flowvec^*))}{\sum_{s\in\rs}y_s(\flowvec^*)G_s(\y(\flowvec^*))}=\frac{\flow_r^*}{\demand},
\end{equation}
which actually coincides with the choice probability definition of Equation (\ref{eqn:AMEVprob2}). Therefore, VI problem (\ref{eqn:suevi}) describes the stochastic user equilibrium with the choice model that is defined by generating function $G$ and generating vector $\y$. 

Such a VI formulation is useful, since the corresponding theory can be used to analyse existence and uniqueness of stochastic user equilibria (as for example done by \citet{Nagurney1998,tr-b14}). In addition, the VI formulation provided the generalized stochastic cost per route -- which should all be equal in equilibrium --, that can be used in assignment algorithms to determine duality gaps. Namely, as $\hat \flowvec\rightarrow\flowvec^*$, then
\begin{equation}
\frac{\sum_{r\in \rs}\hat \flow_r\left(c_r(\hat\flowvec)-\min_{s\in\rs} c_s(\hat\flowvec)\right)}{\sum_{r\in \rs}\hat \flow_r\min_{s\in\rs} c_s(\hat\flowvec)}\rightarrow 0.
\end{equation}  

\subsection{Normalization, identification, and invariance}\label{sec:norm}

The constant in the systematic utility (i.e., the alternative specific constant) of the logit-based (i.e., additive) models has to be normalized because only differences in utility matter. Furthermore, one of the attribute parameters or the scale has to be normalized due to identification. If the cost parameter is normalized to one, all other attribute parameters can  be interpreted as willingness to pay for its attribute. 

The multiplicative models however, do not require the constant in the utility to be normalized. This constant is also multiplied with the error term and thus this term is not equal amongst alternatives. This allows the modeller to use an additional parameter. Similar as for the additive case, one of the attribute parameters or the scale has to be normalized due to identification.

For logit based additive models only differences between utilities matter for the choice probabilities. Similarly, only the ratios between utilities matter for the choice probabilities in the multiplicative models \citep{Fosgerau2009494}. This means that the additive models are invariant under addition with a constant and the multiplicative models are invariant under multiplication with a constant. As the next example shows, it is doubtful that these properties are realistic and they limit both models. For the additive case consider two route sets that only differ by a constant, one with route costs (e.g., travel time) $\left\{1,6\right\}$ and one with $\left\{100,105\right\}$; it is not expected that the choice probabilities are the same for these routes sets. Analogously, route costs $\left\{2,3\right\}$ and $\left\{40,60\right\}$ should neither reflect the same probabilities for the multiplicative case. Nevertheless, since the constant in the systematic utility does not have to be normalized, the property is less restrictive for the multiplicative case. Section \ref{sec:networkchange} provides further analysis.

\section{Multiplicative MEV Models with Explicit Removal of Overlap}

We present an adjustment to the multiplicative model based on the decision rule that travellers only compare the non-overlapping part of routes compared to a reference route. This relaxes the invariance properties to a certain extent, and implies a different choice mechanism for travellers. Probabilities for switching to another route and staying at the reference route are determined. This different decision rule can be retrieved from an alternative utility formulation, which allows a qualitative analysis of the model.  

Using the systematic utility per link, and the multiplicative utility formula, let
\[
U^\Mm_r=\left(\sum_{l\in\LL_r}V_l\right)\rv_r^\Mm, \quad \forall r\in \rs.
\] 
The choice probability for alternative $r\in \rs$ is
%\begin{align}
%P_r^\Mm=&\prob(U_r^\Mm\geq U_p^\Mm, \forall p\neq r)\\=&\prob\left(\left(\sum_{l\in \LL_r}V_l\right)\rv_r^\Mm\geq\left(\sum_{l\in \LL_p}V_l\right)\rv_p^\Mm, \quad \forall p\neq r\right)\\
%=&\prob\left(\sum_{l\in \LL_r\setminus \LL_p}V_l\rv_r^\Mm+\sum_{l\in \LL_r\cap \LL_p}V_l\rv_r^\Mm\geq\sum_{l\in \LL_p\setminus \LL_r}V_l\rv_p^\Mm+\sum_{l\in \LL_r\cap \LL_p}V_l\rv_p^\Mm, \quad \forall p\neq r\right).\label{eqn:componentP}
%\end{align}
\begin{align}
	P_r^\Mm=&\prob(U_r^\Mm\geq U_p^\Mm, \forall p\neq r)\\=&\prob\left(\left(\sum_{l\in \LL_r}V_l\right)\rv_r^\Mm\geq\left(\sum_{l\in \LL_p}V_l\right)\rv_p^\Mm, \quad \forall p\neq r\right)\\
	=&\prob\left(\sum_{l\in \LL_r}V_l\rv_r^\Mm\geq\sum_{i\in \LL_p}V_l\rv_p^\Mm, \quad \forall p\neq r\right)\\
	=&\prob\left(\sum_{l\in \LL_r\setminus \LL_p}V_l\rv_r^\Mm+\sum_{l\in \LL_r\cap \LL_p}V_l\rv_r^\Mm\geq\sum_{l\in \LL_p\setminus \LL_r}V_l\rv_p^\Mm+\sum_{l\in \LL_r\cap \LL_p}V_l\rv_p^\Mm, \quad \forall p\neq r\right).\label{eqn:componentP}
\end{align}
In this formulation, each part of random utility is associated with a part of systematic utility according to the scaling postulate. The systematic utility that is shared amongst a pair of alternatives appears at both sides of the inequality. The shared links for route pair $r,p$ are contained in $\LL_r\cap \LL_p$. An individual will evaluate the links independent of the alternatives they belong to. Therefore we assume that the part of the random utility belonging to these shared links is fully correlated. This means that for all $l \in \LL_r\cap \LL_p$ the terms $V_l\rv_r^\Mm$ and $V_l\rv_p^\Mm$ in Equation (\ref{eqn:componentP}) are equal and can be subtracted from both sides of the inequality. 

This can be derived in an econometrical sound fashion by reconsidering the utility formula based on a reference route. Let $r$ be the reference route, and denote the utility for each route $p \in \rs$ as
\begin{equation}
U_p^{\Md,r}:=\begin{cases}
\sum_{l\in\LL_r}V_l\rv_r^\Mm &\text{if }p=r\\
\sum_{l\in \LL_p\setminus \LL_r}V_l\rv_p^\Mm+\sum_{l\in \LL_r\cap \LL_p}V_l\rv_r^\Mm &\text{otherwise}.
\end{cases}
\end{equation}
The utility that overlaps with the reference route, is multiplied with the error term of the reference route, while the remainder has its own error term. We refer to this type of probabilities as the \emph{\MD case} with reference route $r$. The probability of choosing route $p\in\rs$, conditional on reference route $r$ is  
\begin{equation}
P_p^{\Md,r}=\begin{cases}
\prob\left(\sum_{l\in \Kmin{p}{s}}V_l\rv_p^\Mm
\geq
\sum_{l\in \Kmin{s}{p}}V_l\rv_s^\Mm, \quad \forall s\neq p\right) &\text{if }p=r\\
\prob\left(\sum_{l\in \Kmin{p}{r}}V_l\rv_p^\Mm+\sum_{l\in \LL_r\cap \LL_p}V_l\rv_r^\Mm
\geq
U_s^{\Md,r}, \quad \forall s\neq p\right) &\text{otherwise},
\end{cases}
\end{equation}
with notation $\Kmin{r}{p}=\LL_r\setminus \LL_p$, and the overlapping parts are subtracted for case $p=r$ (see Equation (\ref{eqn:componentP})). First, we analyse the case $p=r$. That choice probability reads that the utility of non-overlapping links of routes $s$ with $p$ is smaller than the utility of the non-overlapping links of route $p$ with each respective $s$. Thus, it is the probability that the traveller does not benefit from swapping some links and change route. Now isolate $\rv_r$ on the left hand side;
\begin{align}
\text{if } p=r \text{, then }P_p^{\Md,r}=&
\prob\left(\left(\sum_{l\in \Kmin{p}{s}}V_l\right)\rv_p^\Mm
\geq
\left(\sum_{l\in \Kmin{s}{p}}V_l\right)\rv_s^\Mm, \quad \forall s\neq p\right)\\
=&\prob\left(\rv_p^\Mm
\leq
\frac{\sum_{l\in \Kmin{s}{p}}V_l}{\sum_{l\in \Kmin{p}{s}}V_l}\rv_s^\Mm, \quad \forall s\neq p\right).\label{eqn:Pi}
\end{align}
With the same algebra that derives Equation (\ref{eqn:MultAdd}) from Equation (\ref{eqn:MultAddStart}), this multiplicative formulation can be transformed to an additive formulation:   
%\begin{align}
%\text{if } p=r \text{, then }P_p^{\Md,r}=&\prob\left(\rv_p^\Mm
%\leq
%\frac{\sum_{l\in \Kmin{s}{p}}V_l}{\sum_{l\in \Kmin{p}{s}}V_l}\rv_s^\Mm, \quad \forall s\neq p\right)\label{eqn:MultiRefRoute}\\
%=&\prob\left(-\ln(\rv_p^\Mm)
%\geq
%-\ln\left(\frac{\sum_{l\in \Kmin{s}{p}}V_l}{\sum_{l\in \Kmin{p}{s}}V_l}\right)-\ln\left(\rv_s^\Mm\right), \quad \forall s\neq p\right)\label{eqn:probMDMEV}
%\end{align}
\begin{align}
	\text{if } p=r \text{, then }P_p^{\Md,r}=&\prob\left(\rv_p^\Mm
	\leq
	\frac{\sum_{l\in \Kmin{s}{p}}V_l}{\sum_{l\in \Kmin{p}{s}}V_l}\rv_s^\Mm, \quad \forall s\neq p\right)\label{eqn:MultiRefRoute}\\
	=&\prob\left(\ln(\rv_p^\Mm)
	\leq
	\ln\left(\frac{\sum_{l\in \Kmin{s}{p}}V_l}{\sum_{l\in \Kmin{p}{s}}V_l}\rv_s^\Mm\right), \quad \forall s\neq p\right)\\
	=&\prob\left(\ln(\rv_p^\Mm)
	\leq
	\ln\left(\frac{\sum_{l\in \Kmin{s}{p}}V_l}{\sum_{l\in \Kmin{p}{s}}V_l}\right)+\ln\left(\rv_s^\Mm\right), \quad \forall s\neq p\right)\\
	=&\prob\left(-\ln(\rv_p^\Mm)
	\geq
	-\ln\left(\frac{\sum_{l\in \Kmin{s}{p}}V_l}{\sum_{l\in \Kmin{p}{s}}V_l}\right)-\ln\left(\rv_s^\Mm\right), \quad \forall s\neq p\right)\label{eqn:probMDMEV}
\end{align}
For the case $p\neq r$ we can show that;  
\begin{align}
\text{if } p\neq r \text{ then }P_p^{\Md,r}=&\prob\left(\sum_{l\in \Kmin{p}{r}}V_l\rv_p^\Mm+\sum_{l\in \LL_r\cap \LL_p}V_l\rv_r^\Mm
\geq
U_s^{\Md,r}, \quad \forall s\neq p\right)\\
=&\prob\left(\frac{\sum_{l\in \Kmin{p}{r}}V_l}{\sum_{l\in \Kmin{r}{p}}V_l}\rv_p^\Mm\leq\rv_r^\Mm \wedge \left(\frac{\sum_{l\in \Kmin{p}{s}}V_l}{\sum_{l\in \Kmin{s}{p}}V_l}\rv_p^\Mm\leq\frac{\sum_{l\in \Kmin{s}{p}}V_l}{\sum_{l\in \Kmin{p}{s}}V_l}\rv_s^\Mm,\forall s\neq r,p\right)
\right)\label{eqn:nonTriv}
\end{align}
This equality is not trivial, one should check that Equation (\ref{eqn:nonTriv}) holds or fails for all six orderings of $U_r,U_p$, and $U_s$ by using and/or substituting the condition in Equation (\ref{eqn:MultiRefRoute}). So, if route $p$ is chosen under reference route $r$, then it is beneficial to switch from route $r$ to route $p$, and this improvement is larger than to switch to any other route $s$. Equation (\ref{eqn:nonTriv}) can be rewritten in an additive form similar to (\ref{eqn:probMDMEV}) by taking the logarithm transformation. Assume -- similar as in the \M models -- that in the derived additive forms $\left(-\ln(\rv_1^\Mm), \ldots, -\ln(\rv_{|\rs|}^\Mm)\right)$ follows \MEV\ distribution $\RV^\MEV$, then we have another type of \GMEV\ route choice models. The generating vector is reference route specific, denoted with $\yMDr$, and given by
\begin{equation}
\yyMDr_{p}=
\begin{cases}
1 &\text{if }r=p\\
\frac{\sum_{l\in \Kmin{r}{p}}V_l}{\sum_{l\in \Kmin{p}{r}}V_l} &\text{otherwise} 
\end{cases}
, \quad \forall r,p\in\rs,
\end{equation}
where we used that $e^0=1$ and \[e^{-\ln\left(\frac{\sum_{l\in \Kmin{p}{r}}V_l}{\sum_{l\in \Kmin{r}{p}}V_l}\right)}=\frac{\sum_{l\in \Kmin{r}{p}}V_l}{\sum_{l\in \Kmin{p}{r}}V_l}.\]
Given any generating function $G$, reference route $r$ and \emph{\MD generating vector} $\yMDr$ the choice probabilities can be derived with Equation (\ref{eqn:Gprobs}). By then applying Euler's homogeneous function theorem, the probability of choosing route $p\in\rs$ with reference route $r\in \rs$ is:
\begin{equation}
P_p^{\Md,r}=\frac{\yyMDr_p G_p(\yMDr)}{G_r(\yMDr)+\sum_{\{p\in \rs|p\neq r\}}\frac{\sum_{l\in \Kmin{p}{r}}V_l}{\sum_{l\in \Kmin{r}{p}}V_l}G_p(\yMDr)}.\label{eqn:MDprobConditional}
\end{equation}
Note that if no overlap exists the \MD models collapse to the \M models\ (i.e., they return the same probabilities). 

The previous analysis was based on one reference route, but for applications there is not always a (single) reference route available since the (current) reference of the travellers is not known. The final step is to handle this uncertainty. Denote $\Pref{r}$ as the probability that route $r$ is the reference route. The final choice probability for route $p$ then becomes 
\begin{equation}
P_p^\Md=\sum_{r\in\rs} P_p^{\Md,r}\Pref{r}.\label{eqn:MDprob}
\end{equation} 
We consider three ways to determine the $\Pref{r}$. First, it is possible to assign every route as reference route with the same probability (i.e., $\Pref{r}=1/|\rs|,\forall r\in \rs$). This seems only realistic when no irrelevant routes exist in the choice set, since it is not plausible that travellers use an irrelevant route as reference. Second, it is very natural to set the probability that a route is chosen equal to the probability that a route is the reference route. Then Equation (\ref{eqn:MDprob}) becomes a system of equations: $P_p^\Md=\sum_{r\in\rs} P_p^{\Md,r}P_p^\Md$ (where $P_p^\Md,p\in\rs$ are the unknowns). Furthermore, $\sum_{p\in\rs}P_p^\Md=1$, thus this system of equations can be identified as a Markov chain. Since $P_p^{\Md,r}>0, \forall p,r\in\rs$, a steady state exists which can be found by solving the system of equations. A third possibility is to fix one route (e.g., the fastest in free-flow conditions) as the reference route $r$, and to only determine the choice probabilities of $\Md,r$. This means that $\Pref{r}=1$ for exactly one $r\in\rs$. This might not be the realistic for standard equilibrium models since this creates asymmetry. However, such an approach would be very feasible for en-route decisions, where the current route serves as reference point. Also in day-to-day models, the previously chosen trip can be a natural reference route.

This new family of \Md\ does not account automatically for all overlap with the multinomial generating function. Any network with overlap in which all routes have the same length, will give equal choice probabilities; the network structure does not have any influence on the choice probabilities. Also, the conditional choice situation with reference route $r$ does account for overlap between $r$ and all other alternatives, but cannot capture dependencies between all these other alternatives. So, a specific generating function that can handle these dependencies has to be specified per reference route. 

\subsection{Model instances}

Using the four available generating functions we can derive several models based on multiplicative utility formulas based on reference routes. Here we only present the choice probabilities of route $r$ conditional on $r$ being the reference route, i.e., using $p=r$ in Equation (\ref{eqn:MDprobConditional})\footnote{Under the Markov chain assumption, these are the probabilities for staying in the same state.}. The derived choice probabilities seem to be rather complex; however, that is only due to the asymmetry of the generating vector. Basically, the choice probabilities of a $\Md,r$-model are not more complex than those for a \A\ or \M model, and they can be derived easily for applications.   
\begin{description}
\item[\MD\MNL] The disadvantage of independence in basic MNL is inherited. Not all overlap is explicitly removed in the \MD case, dependencies between non-reference routes maintain. The choice probabilities are
\begin{equation}
P_r^{\MD\MNL,r}(G^\MNL;\yMDr)=\frac{1}{1+\sum_{p\in\rs\setminus\{r\}}\left(\frac{\sum_{l\in \Kmin{r}{p}}V_l}{\sum_{l\in \Kmin{p}{r}}V_l}\right)^\scale}
,\quad \forall r\in \rs.
\end{equation}
\item[\MD\PSL]
The path-size choice probabilities are
\begin{equation}
P_r^{\MD\PSL,r}(G^\PSL;\yMDr)=\frac{\ps_r^\pspar}{\ps_r^\pspar+\sum_{p\in\rs\setminus\{r\}}\ps_p^\pspar\left(\frac{\sum_{l\in \Kmin{r}{p}}V_l}{\sum_{l\in \Kmin{p}{r}}V_l}\right)^\scale}
,\quad \forall r\in \rs.
\end{equation}
As it is possible to convert the path-size factor back into the utility formulation for the \A\PSL\ and \M\PSL\ models \citep{Kitthamkesorn2013}, this is also possible for the \MD\PSL\ model. This will lead to: 
\begin{equation}
U_p^{\MD\PSL,r}:=\begin{cases}
\frac{\sum_{l\in\LL_r}V_l\rv_r^\Mm}{\ps_r^\pspar} &\text{if }p=r\\
\frac{\sum_{l\in \LL_p\setminus \LL_r}V_l\rv_p^\Mm}{\ps_p^\pspar}+\frac{\sum_{l\in \LL_r\cap \LL_p}V_l\rv_r^\Mm}{\ps_r^\pspar} &\text{otherwise}.
\end{cases}
\end{equation}
Having a reference route allows one to revisit the used path-size factors. These can be obtained by excluding all links of the reference route. This will lead to a  formulation with reference route specific path-size parameters. Similar to equation (\ref{eqn:pathsize}), the formulation for the reference route $r$ specific path-size factor $\ps_{p,r}$ for route $p$ is
\begin{equation}\label{eqn:pathsizeMD}
\ps_{p,r}=\begin{cases}\frac{\sum_{l\in \Kmin{p}{r}}\frac{V_l}{\#_l}}{\sum_{l\in \Kmin{p}{r}}V_l}, & \text{if } p\neq r\\
1&\text{if } p=r.\end{cases}
\end{equation}

\item[\MD\PCL]
The paired combinatorial choice probabilities for in the \MD case are
\begin{align}
& P_r^{\MD\PCL,r}(G^\PCL;\yMDr)=\nonumber\\
& \frac{
 \sum_{p\in\rs\setminus\{r\}}  \left(1+\left(\frac{\sum_{l\in \Kmin{r}{p}}V_l}{\sum_{l\in \Kmin{p}{r}}V_l}\right)^\frac{\scale}{1 - \cf_{rp}}\right)^{-\cf_{rp}}
}{
 \sum_{r'\in\rs\setminus\{r\}}
 \left(1+\left(\frac{\sum_{l\in \Kmin{r}{r'}}V_l}{\sum_{l\in \Kmin{r'}{r}}V_l}\right)^\frac{\scale}{1 - \cf_{rr'}}\right)^{1-\cf_{rr'}}+
 \sum_{p\in\rs\setminus\{r,r'\}}\left(\left(\frac{\sum_{l\in \Kmin{r}{r'}}V_l}{\sum_{l\in \Kmin{r'}{r}}V_l}\right)^\frac{\scale}{1 - \cf_{r'p}}
 +\left(\frac{\sum_{l\in \Kmin{r}{p}}V_l}{\sum_{l\in \Kmin{p}{r}}V_l}\right)^\frac{\scale}{1 - \cf_{r'p}}\right)^{1-\cf_{r'p}}
}
,\quad \forall r\in \rs.
\end{align}
This formula looks complicated, as the cumbersome indexing is required because the non-symmetric generating vector $\yMDr$ is substituted, and then some terms collapse to 1. Without this substitution the probabilities are
\begin{equation}
 P_r^{\MD\PCL,r}(G^\PCL;\yMDr)=
 \frac{
 \sum_{p\in\rs\setminus\{r\}}  \left({\yyMD_{rr}}^\frac{\scale}{1 - \cf_{rp}}+{\yyMD_{rp}}^\frac{\scale}{1 - \cf_{rp}}\right)^{-\cf_{rp}}
}{
 \sum_{r'\in\rs}\sum_{p\in\rs\setminus\{r'\}}\left({\yyMD_{rr'}}^\frac{\scale}{1 - \cf_{r'p}}+{\yyMD_{rp}}^\frac{\scale}{1 - \cf_{r'p}}\right)^{1-\cf_{r'p}}
}
,\quad \forall r\in \rs.
\end{equation}
\item[\MD\LNL]
The link-nested choice probabilities in the \MD case  are
\begin{equation}
P_r^{\MD\LNL,r}(G^\LNL;\yMDr)=
\frac{
 \sum_{m\in\LL} \ic_{mr} 
 \left(\ic_{mr}+\sum_{p\in\rs\setminus\{r\}}\ic_{mp}\left(\frac{\sum_{l\in \Kmin{r}{p}}V_l}{\sum_{l\in \Kmin{p}{r}}V_l}\right)^{\scale_m}\right)^{\frac{\scale}{\scale_m}-1}
}{
 \sum_{m\in\LL} \left(\ic_{mr}+\sum_{p\in\rs\setminus\{r\}}\ic_{mp}\left(\frac{\sum_{l\in \Kmin{r}{p}}V_l}{\sum_{l\in \Kmin{p}{r}}V_l}\right)^{\scale_m}\right)^\frac{\scale}{\scale_m}
}
,\quad \forall r\in \rs.
\end{equation}
\end{description}

\subsection{Model properties}

Since for the \Md\ case utilities are sums of two variates of the MEV distributions, it is not possible to obtain a convenient closed-form for the joint density function of the utilities. Therefore, we cannot derive the expected maximum utility. However, from the marginal distributions of $\RV^\Mm$ it is possible to obtain the expected value and variance of the \Md\ utilities conditional on a reference route, which are:   

\begin{equation}
\E(U_p^{\Md,r})=\begin{cases}
\frac{\sum_{l\in\LL_r}V_l}{G(\one_r)^{1/\scale}}\Gamma\left(1+\frac{1}{\scale}\right) &\text{if }p=r\\
\left(\frac{\sum_{l\in \LL_p\setminus \LL_r}V_l}{G(\one_p)^{1/\scale}}+\frac{\sum_{l\in \LL_r\cap \LL_p}V_l}{G(\one_r)^{1/\scale}}\right)\Gamma\left(1+\frac{1}{\scale}\right) &\text{otherwise}
\end{cases}, \text{ and,}
\end{equation}

\begin{equation}
\Var(U_p^{\Md,r})=\begin{cases}
\frac{\left(\sum_{l\in\LL_r}V_l\right)^2}{G(\one_r)^{2/\scale}}\left(\Gamma\left(1+\frac{2}{\scale}\right)-\Gamma\left(1+\frac{1}{\scale}\right)^2\right) &\text{if }p=r\\
\left(\frac{\left(\sum_{l\in \LL_p\setminus \LL_r}V_l\right)^2}{G(\one_p)^{2/\scale}}+\frac{\left(\sum_{l\in \LL_r\cap \LL_p}V_l\right)^2}{G(\one_r)^{2/\scale}}\right)\left(\Gamma\left(1+\frac{2}{\scale}\right)-\Gamma\left(1+\frac{1}{\scale}\right)^2\right) &\text{otherwise}.
\end{cases}
\end{equation}

Furthermore, when the reference route is fixed (i.e., $\Pref{r}=1$ for exactly one $r\in\rs$), then SUE formulation (\ref{eqn:suevi}) also holds for the $\Md,r$-case.

\subsection{Simple network}

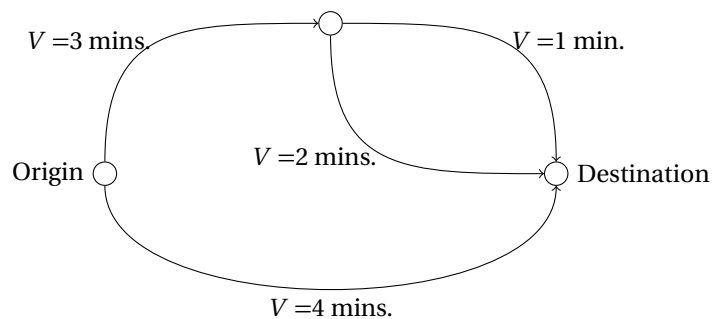
\begin{figure}[htbp]
\centering
\begin{tikzpicture}
\node[label=left:Origin,circle,draw]
(n1) at (0,0) {};
\node[circle,draw]
(n2) at (3,2) {};
\node[circle,draw,label=right:Destination]
(n3) at (6,0) {};

\draw[->] (n1) ..controls +(90:2cm) and +(180:2cm)..  node[above,left] {$V=$3 mins.}  (n2);
\draw[->] (n2) ..controls +(0:2cm) and +(90:2cm)..  node[above,right] {$V=$1 min.}  (n3);
\draw[->] (n2) ..controls +(-90:2cm) and +(180:2cm)..   node[below,left] {$V=$2 mins.} (n3);
\draw[->] (n1) ..controls +(-90:2cm) and +(-90:2cm)..  node[below] {$V=$4 mins.} (n3);
\end{tikzpicture}
\caption{Simple three route overlap network.}
\label{fig:simplenetwork}
\end{figure}

To provide more insight in the working of \MD models, the \MD\MNL\ probabilities are provided for a simple network. Figure \ref{fig:simplenetwork} shows an origin and destination with three routes in between. Systematic utility is assumed to equal foreseen travel time, and no route specific constant is included. The upper and middle routes have overlap, and the upper route is faster. The lower route has no overlap with the other two, and is equally fast as the upper route. Table \ref{tab:simplenetwork} shows the choice probabilities for each reference route based on $\scale=1$, as well as, the solutions based on equal reference route probabilities and the Markov chain approach.  

\begin{table}
\caption{Choice probabilities for the \MD\MNL\ model on the simple network. The probabilities conditional on reference routes, as well as the two solution methods are provided.}
	\centering
\begin{tabular}{l|lll}
\toprule  & Upper & Middle & Lower \\ 
\midrule Ref. route = Upper & $\frac{2}{5}$ & $\frac{1}{5}$ & $\frac{2}{5}$ \\[5pt] 
Ref. route = Middle & $\frac{8}{17}$ & $\frac{4}{17}$ & $\frac{5}{17}$ \\[5pt]
Ref. route = Lower & $\frac{5}{14}$ & $\frac{4}{14}$ & $\frac{5}{14}$ \\[5pt]
Equal $\Pref\cdot$ solution & $\approx 0.409$ & $\approx 0.240$ & $\approx 0.350$ \\[5pt] 
Markov chain solution & $\approx 0.401$ & $\approx 0.239$ & $\approx 0.359$ \\ 
\bottomrule 
\end{tabular}  
\label{tab:simplenetwork}
\end{table}

As expected, the slowest middle route has lowest choice probability. This choice probability is significantly smaller than it would be under \A\ and \M models, which is realistic since changing from the upper to the middle route means that the non-overlapping travel time doubles. Also, the two different solution methods for dealing with multiple reference routes do not differ that much. 

On the other hand, there is a higher preference for the upper route than for the lower route while they have equal travel times. Based on the theory in Section \ref{sec:rum}, one would argue that the opposite should be true. Having the upper or lower route as the reference route, causes no difference in choice probabilities between the upper and lower route. Only because the middle route exists, and since it is more 'profitable' (i.e. the factor 2 between the non-overlapping parts between middle and lower is higher than the factor $4/3$  between middle and lower route) to switch to the upper route, there is a slight preference for the upper route. This means that by having routes included as reference routes, the final probabilities change (similar to what \citet{bliemertrr2008} show). 

So, this simple example provides insight in the working of the \MD models, but does not show its full potential. Clearly, the \MD models have different behaviour than all other known route choice models.  The next two sections provide more insight in the differences between the \A, \M and \MD models -- and their suitability for traffic assignment applications --. First, we analyse their basic behaviour under simple changes in the network configuration, and second, we analyse all models under more complex network variations. The latter quantitative test shows which model can best approximate the generic utility formulation of Section \ref{sec:rum}. These Sections provide the full potential of the \MD models. 

\section{Basic model behaviour under simple network changes}\label{sec:networkchange}

In this section we discuss the change in choice behaviour under three simple network adjustments. Consider the four networks depicted in the first column of Table \ref{tab:networkvariations}. Network A is the basis and has 2 routes with different lengths of which the final parts overlap. In network B the overlapping part is extended, while in networks C1 and C2 the non-overlapping parts are extended by respectively adding a constant length and by multiplying their lengths. The models are so simple that, regardless the \emph{exact} choice behaviour, the trend in choice probabilities for the two routes is known. The trend can either be that the probabilities remain equal, converge, or diverge. It is desired that choice models can reproduce the expected trend for each network change. However, not all trends can be reproduced by all choice models. Therefore, we analyse the behaviour of the three basic (-\MNL) models. We take advantage of having only two routes for the \MD\MNL\ model here, which avoids dealing with reference routes. However, similar results are obtained when extending to more than two routes, as is analysed in the next section.

  \begin{table}
  \caption{\textbf{Behaviour under changing networks.} Networks B, C1 and C2 are slight variations on network A. Each network has two routes; the table shows the trend of the route probabilities when one switches for network A to any of the other three networks. They can either converge, diverge, or remain equal (depicted with arrows). The expected trend and achievable trends for three models are provided (see main text).}
  	\centering
  	 
  \begin{tabular}{ll|m{2cm}m{1.5cm}m{1.5cm}m{1.5cm}}
  \toprule  A & \parbox{5cm}{\begin{tikzpicture}
  	\node[circle,draw]
  	(n1) at (0,0) {};
  	\node[circle,draw,scale=0.5]
  	(n2) at (1.5,0) {};
  	\node[circle,draw]
  	(n3) at (3,0) {};
  	
  	\draw[->] (n1) ..controls +(0:1cm) and +(180:1cm)..  node[above] {1}  (n2);
  	\draw[->] (n1) ..controls +(-45:1cm) and +(225:1cm)..  node[below] {2}  (n2);
  	\draw[->] (n2) -- node[above] {1} (n3);
  	\end{tikzpicture}}  & Expected choice behaviour & \A\MNL & \M\MNL & \MD\MNL \\ 
  \midrule B & \parbox{5cm}{\begin{tikzpicture}
  \node[circle,draw]
  (n1) at (0,0) {};
  \node[circle,draw,scale=0.5]
  (n2) at (1.5,0) {};
  \node[circle,draw]
  (n3) at (4.5,0) {};
  
  \draw[->] (n1) ..controls +(0:1cm) and +(180:1cm)..  node[above] {1}  (n2);
  \draw[->] (n1) ..controls +(-45:1cm) and +(225:1cm)..  node[below] {2}  (n2);
  \draw[->] (n2) -- node[above] {2} (n3);
  \end{tikzpicture}} & $=$ & $=$ \quad\cmark & $\converge$ \quad\xmark & $=$ \quad\cmark \\ 
  \cmidrule(lr){1-2} C1 & \parbox{5cm}{\begin{tikzpicture}
  \node[circle,draw]
  (n1) at (0,0) {};
  \node[circle,draw,scale=0.5]
  (n2) at (3,0) {};
  \node[circle,draw]
  (n3) at (4.5,0) {};
  
  \draw[->] (n1) ..controls +(0:1cm) and +(180:1cm)..  node[above] {2}  (n2);
  \draw[->] (n1) ..controls +(-45:1cm) and +(225:1cm)..  node[below] {3}  (n2);
  \draw[->] (n2) -- node[above] {1} (n3);
  \end{tikzpicture}} & \converge  & $=$ \quad\xmark & \converge \quad\cmark & \converge \quad\cmark \\ 
  \cmidrule(lr){1-2} C2 & \parbox{5cm}{\begin{tikzpicture}
  \node[circle,draw]
  (n1) at (0,0) {};
  \node[circle,draw,scale=0.5]
  (n2) at (3,0) {};
  \node[circle,draw]
  (n3) at (4.5,0) {};
  
  \draw[->] (n1) ..controls +(0:1cm) and +(180:1cm)..  node[above] {2}  (n2);
  \draw[->] (n1) ..controls +(-60:1.5cm) and +(240:1.5cm)..  node[below] {4}  (n2);
  \draw[->] (n2) -- node[above] {1} (n3);
  \end{tikzpicture}} & \diverge & \diverge \quad\cmark & \diverge \quad\cmark & \diverge \quad\cmark \\ 
  \bottomrule 
  \end{tabular}\\
   
  \label{tab:networkvariations}
  \end{table}
  
Consider the switch from network A to network B. Since the non-overlapping parts remain equal, the choice probabilities for both routes obviously also remain equal. In the \A\MNL\ model only the difference between the routes matters, and since this difference does not change, the probabilities also remain equal for that model. For the \M\MNL\ model however, only the ratio between the routes matters, and this ratio changes. Because the non-overlapping part increases, the (lower:upper)-ratio  decreases, and the choice probabilities will converge. The expected behaviour can not be reproduced by any \M\MNL\ model instance. The \MD\MNL\ model in the end, merely considers the ratio of the difference between routes, i.e., the ratio of the non-overlapping parts. Since the non-overlapping part does not change, the choice probabilities do not change. This holds for all \MD\MNL\ model instances 

Consider the switch from network A to network C1. When a constant length is added to both routes, they become more similar the choice probability of the shortest route increases. So, it is expected that the route probabilities will converge\footnote{As an easy example, consider routes with length $x$ and $x+1$; the choice probabilities will become 50 \%-50\% if $x\rightarrow\infty$ and 100\%-0\% if $x\rightarrow 0$}. The difference between the two routes will not change, therefore no \A\MNL\ model instance can reproduce the expected behaviour. On the other hand, the \M\MNL\ and \MD\MNL\ models will reproduce the expected behaviour since the (lower:upper)-ratio  between routes (in- and excluding the overlapping part) decreases. The probabilities converge for all \M\MNL\ and \MD\MNL\ model instances.

Consider the switch from network A to network C2. When a both routes are multiplied by the same factor, the absolute detour of the longest route will increase, and will therefore be chosen less. So, it is expected that route probabilities will diverge\footnote{As an easy example, consider routes with length $x$ and $2\times x$; the choice probabilities will become 50 \%-50\% if $x\rightarrow 0$ and 100\%-0\% if $x\rightarrow \infty$}. The difference between the two routes will increase, therefore the probabilities in the \A\MNL\ models will diverge, which coincides with the expected behaviour. The ratio of the non-overlapping parts of the routes remains equal, but due to the overlapping part, the ratio between the whole routes changes. This increase in the (lower:upper)-ratio will lead to the desired diverging probabilities in the \M\MNL\ model. Of course, the divergence `speed' depends on the length of the overlapping part. However, this can be adjusted by using the constant in the systematic utility specification. This constant is also the reason that the diverging probabilities are obtainable in the \MD\MNL\ model; here it also completely determines the divergence `speed', which is an advantage compared to \M\MNL. All \M\MNL\ model instances diverge for this network example, but the behaviour remains dependent on the overlapping length. Almost all \MD\MNL\ model instances will have the expected behaviour, except for those with a constant equal to zero.     

This analysis, summarized in Table \ref{tab:networkvariations}, on the most simple and basic network variation shows that only the \MD\MNL\ model can reproduce all expected behaviours. We believe no simple network change exists of which the expected behaviour can be captured by \A\MNL\ or \M\MNL, but not by \MD\MNL. Since real networks do not consist of merely simple network changes, the next section analyses the competitiveness of all models on more comprehensive networks based on the route utility formulation of section \ref{sec:rum}, including the \PSL, \PCL, and \LNL\ variants, and including more than two routes which requires the use of reference routes for the \MD models.

\section{Network Example}\label{sec:network}

The route choice models are applied on a network with three routes to show the advantage of the \MD models compared to the others. \citet{bliemertrr2008} compare route choice models on different route sets. They demonstrate that existing closed-form route choice models (except for MNL) are sensitive to irrelevant route alternatives. Hence, if the route sets in application are different from the route sets in estimation, the results may be poor. This network example considers very different choice situations and the models are estimated and validated on different subsets of the data. 

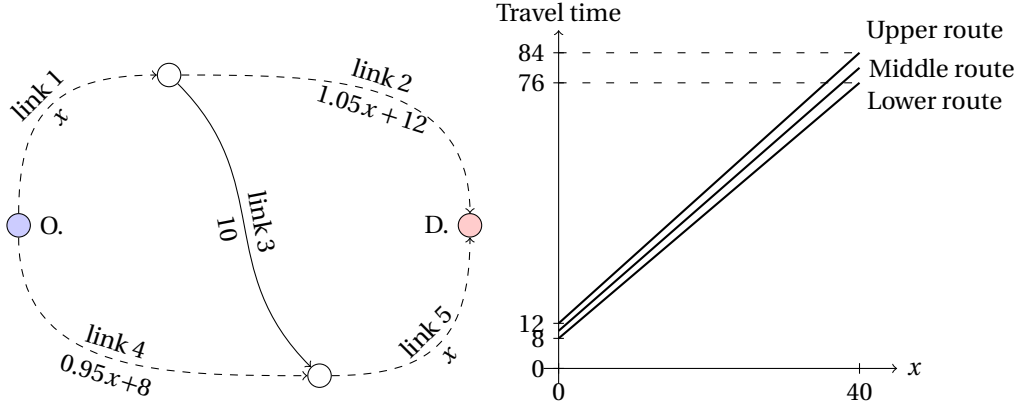
\begin{figure}[bhtp]
\centering
\begin{tikzpicture}
\node[label=right:O.,circle,draw,fill=blue!20]
(n1) at (0,0) {};
\node[circle,draw]
(n2) at (2,2) {};
\node[circle,draw]
(n3) at (4,-2) {};
\node[circle,draw,label=left:D.,fill=red!20]
(n4) at (6,0) {};

\draw[->,dashed] (n1) ..controls +(90:1.5cm) and +(180:1.5cm)..  node[above,sloped] {link 1}  node[below,sloped] {$x$} (n2);
\draw[->,dashed] (n2) ..controls +(0:2cm) and +(90:2cm)..  node[above,sloped] {link 2} node[below,sloped] {$1.05x+12$} (n4);
\draw[->,dashed] (n1) ..controls +(-90:2cm) and +(180:2cm)..  node[above,sloped] {link 4} node[below,sloped] {$0.95x$+8} (n3);
\draw[->] (n2) ..controls +(-45:2cm) and +(135:2cm)..  node[above,sloped] {link 3} node[below,sloped] {10} (n3);
\draw[->,dashed] (n3) ..controls +(0:1.5cm) and +(-90:1.5cm)..  node[above,sloped] {link 5} node[below,sloped] {$x$} (n4);

%\coordinate
%(n4) at (9,1.5) ;
%\coordinate
%(n5) at (10,2.5) ;
%\coordinate
%(n6) at (12,1.5) node[right of=n6,right=-1cm,align=left] {Upper route\\$x+10$};
%
%\draw[dashed] (n4) ..controls +(90:.75cm) and +(180:.75cm)..  (n5);
%\draw[->] (n5) ..controls +(0:1cm) and +(90:1cm)..  (n6);
%\draw[->] (n5) ..controls +(-90:1cm) and +(180:1cm)..  (n6);
%\draw[->,dashed] (n4) ..controls +(-90:1cm) and +(-90:1cm)..  (n6);

%\coordinate
%(n7) at (9,0) ;
%\coordinate
%(n81) at (10,1) ;
%\coordinate
%(n82) at (11,-1) ;
%\coordinate
%(n9) at (12,0) node[right of=n9,right=-1cm,align=left] {Middle route\\$2x+5$};
%
%\draw[dashed] (n7) ..controls +(90:.75cm) and +(180:.75cm)..  (n81);
%\draw (n81) ..controls +(-45:1cm) and +(135:1cm)..  (n82);
%%\draw[->] (n8) ..controls +(0:1cm) and +(90:1cm)..  (n9);
%\draw[->,dashed] (n82) ..controls +(0:.75cm) and +(-90:.75cm)..  (n9);
%%\draw[->,dashed] (n7) ..controls +(-90:1cm) and +(-90:1cm)..  (n9);
%
%\coordinate
%(n10) at (9,-1.5) ;
%\coordinate
%(n11) at (11,-2.5) ;
%\coordinate
%(n12) at (12,-1.5) node[right of=n12,right=-1cm,align=left] {Lower route\\$x+20$} ;
%
%%\draw[dashed] (n10) ..controls +(90:1cm) and +(180:1cm)..  (n11);
%%\draw[->] (n11) ..controls +(0:1cm) and +(90:1cm)..  (n12);
%%\draw[->] (n11) ..controls +(-90:1cm) and +(180:1cm)..  (n12);
%\draw (n10) ..controls +(-90:1cm) and +(180:1cm)..  (n11);
%\draw[->,dashed] (n11) ..controls +(0:.75cm) and +(-90:.75cm)..  (n12);

\end{tikzpicture}
\begin{tikzpicture}[scale=1]

  \draw[->] (-0.2,0) -- (4.5,0) node[right] {$x$};
  \draw[->] (0,-0.2) -- (0,4.5) node[above] {Travel time};

	\draw[shift={(0,0)}] (0pt,2pt) -- (0pt,-2pt) node[below] {$0$};
	\draw[shift={(4,0)}] (0pt,2pt) -- (0pt,-2pt) node[below] {$40$};
  \foreach \y/\ytext in {0/0,0.4/8,0.6/12,3.8/76,4.2/84}
    \draw[shift={(0,\y)}] (2pt,0pt) -- (-2pt,0pt) node[left] {$\ytext$};
	\draw[loosely dashed,thin] (0,4.2) -- (4,4.2);
	\draw[loosely dashed,thin] (0,3.8) -- (4,3.8);
  \draw[thick] (0,0.4) --  (4,3.8) node[right=1cm,below] {Lower route};
  \draw[thick] (0,0.6) --  (4,4.2) node[right=1cm,above] {Upper route};  
  \draw[thick] (0,0.5) --  (4,4) node[right=0cm] {Middle route};   
\end{tikzpicture}
\caption{Network and travel times. Dashed links have variable costs. }
\label{fig:network}
\end{figure}

In Figure \ref{fig:network} the network with link travel times and routes is depicted. This network has been carefully constructed to put the choice models under stress, namely with and without overlap, and with short and long distance OD pairs. Also, the relative influence of the analyst and random foreseen travel time changes. The travel times of links 1, 2, 4 and 5 depend on variable $x$ and by varying $x$ different choice situations are created. Each value of $x$ can represent a different OD pair in a transport network. The graph shows how route travel times increase with $x$; for $x=0$ the difference (respectively ratio) between the slowest and fastest routes is four minutes (respectively 0.67), and this increases to eight minutes (respectively 0.90) for $x=40$. Furthermore, overlap occurs on links 1 and 5.  

For the experiment, the probabilities from the MNP model -- that can properly handle both the analyst error and random foreseen travel time, and the desired variance-covariance structure based on the scaling postulate -- are the ground truth. Thus, the route utilities are jointly distributed following a multivariate normal distribution. This distribution is specified in line with the random route utility of Section \ref{sec:rum}. Assume that there is no systematic utility other than travel time, then identification leads to $V^0_r=0$ for the \A models, and $V^0_r=c$ for the \M\ and \MD models (see Section \ref{sec:norm}). For the normalization set the expected analyst error $\E(\rv_r)$ to $0$ and the travel time parameter $\beta$ to $-1$. For the proportionality parameter (regarding the standard deviation of the foreseen travel time), linear regression on the OViN data of Figure \ref{fig:TTscatter} leads to  $\theta=0.3859$, and linear regression on the route survey leads to  $\theta=0.1301$; however, the first value is too high (see Section \ref{sec:rum:lin}) and the latter too low (since it is the response of only one traveller), thus we assume $\theta=0.2$. The covariances between the random foreseen travel times are assumed to be based on the arithmetic mean, see Equation (\ref{eqn:covAm}). Note that these definitions of the (co)variances are different from those in the literature that assume proportionality between variance and mean. The standard deviation of the analyst error is set to $\stdev(\rv_r)=10$ (minutes). This leads to the following multivariate normal utility distribution: 
\begin{equation}
\mathcal{N}\left(\underbrace{-\left(\begin{array}{c}
2.05x+12 \\ 
2x+10 \\ 
1.95x+8
\end{array}\right)}_{\E(\boldsymbol{U})} ,\underbrace{0.2^2\left(\begin{array}{ccc}
(2.05x+12)^2 & 2.025x^2+11x & 0 \\ 
2.025x^2+11x & (2x+10)^2 & 1.975x^2+9x \\ 
0 & 1.975x^2+9x & (1.95x+8)^2
\end{array} \right)}_{\text{covariance matrix of }\boldsymbol{\TT}}+\underbrace{\text{diag}\left(\begin{array}{c}
100 \\ 
100 \\ 
100
\end{array}\right)}_{\text{variance matrix of }\boldsymbol{\rv}}\right).\label{eqn:MVN}
\end{equation}  

While $x$ increases, not only the distance between origin and destination increases, but also the influence of the randomness from foreseen travel time compared to the randomness from the analyst error increases. For $x=0$ only 3.8 percent of the variance of utility is due to the foreseen travel time; however, for $x=40$, 76.4 percent of the variance is due to the foreseen travel time. Furthermore, the part of utility with overlap increases with $x$. In the development of this network example, the link travel times were chosen given the assumptions on the errors, and such that the route choice probabilities are more or less stable.  

\begin{figure}[htbp]
\centering
\includegraphics[width=7.5cm]{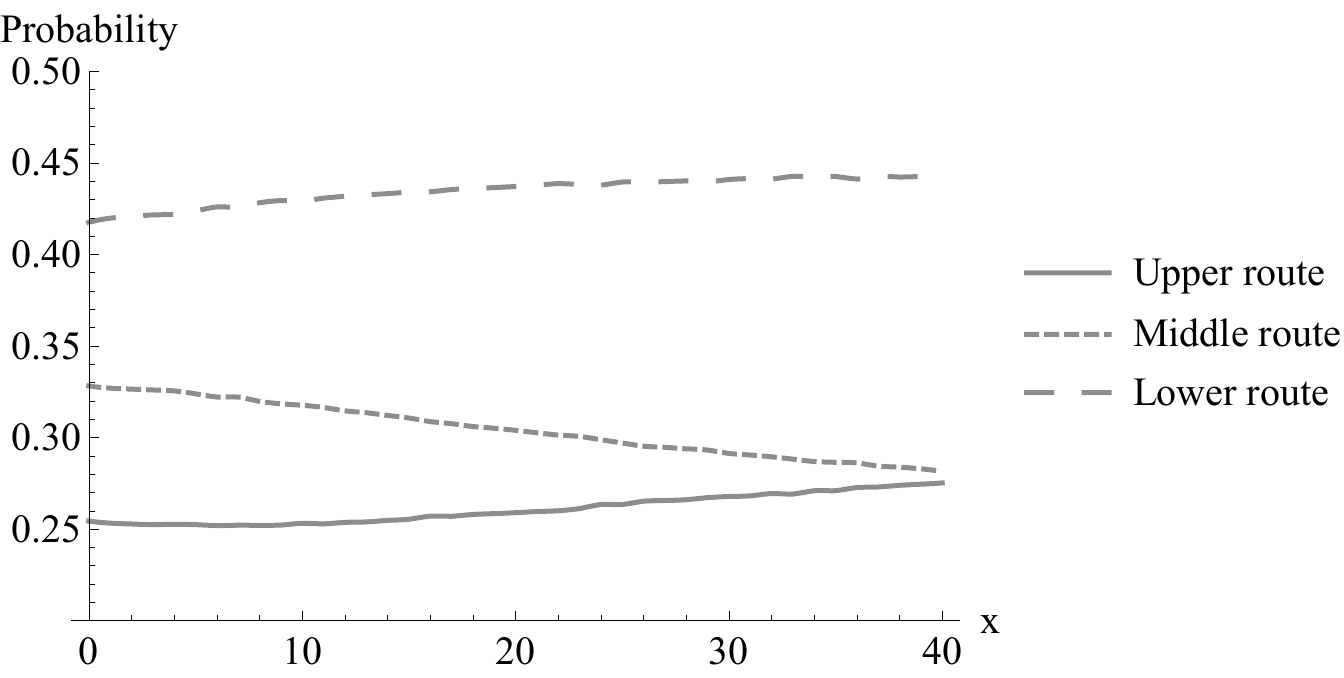}
\caption{Ground truth probabilities based on a million multivariate normal samples for each $x\in\{0,\ldots,40\}$.}
\label{fig:groundtruth}
\end{figure}

For all $x\in \{0,1,\ldots,40\}$ the choice situation is simulated a million times by sampling from the multivariate normal distribution (Equation (\ref{eqn:MVN})), see Figure \ref{fig:groundtruth}. The models are estimated twice; once on the dataset $x=\{5,\ldots,15\}$, and once on dataset $x=\{25,\dots,35\}$ using log-likelihood maximization. So, eleven million `observations' are used to estimate each model. For the \MD\PSL\ the path-size formulation of equation (\ref{eqn:pathsize}) is chosen\footnote{This is chosen since equation (\ref{eqn:pathsizeMD}) leads to much more path-size factors}. For the \LNL\ models link specific scales are estimated for links 1 and 5; the other links -- read nests -- contain only one route and thus `collapse': $(y^{\scale_l})^{\scale/\scale_l}=y^\scale$. The Markov chain approach for reference routes is chosen for \MD\MNL\ and \MD\PSL, and the `equal probability reference route' approach is chosen for \MD\PCL\ and \MD\LNL\footnote{For these two models the software (Wolfram Mathematica) couldn't find analytical soltions of the choice probabilities using the Markov chain with unknown $x$. General applications do not have this parametrization on $x$, and thus shouldn't be problematic.}.  For the validation, the models' log-likelihood on the other dataset is determined, so the parameters from the estimation on $x=\{5,\ldots,15\}$ are applied to $x=\{25,\dots,35\}$, and vice versa. 

\begin{table}
\caption{Parameter estimates for every model for the two datasets}
	\centering
		\begin{tabular}{l|llp{2.0cm}|llp{2.0cm}}
		\toprule
		Model & \multicolumn{6}{c}{Parameter estimates} \\
		 & \multicolumn{3}{c|}{Dataset $x=\{5,\ldots,15\}$} & \multicolumn{3}{c}{Dataset $x=\{25,\ldots,35\}$}\\
		& Scale $\scale\approx$ & Constant $\constant\approx$ & Other & Scale $\scale\approx$ & Constant $\constant\approx$ & Other\\
		\cmidrule(r){1-7}
		\A\MNL & 0.107 & & &  0.0699 & & \\
		\A\PSL & 0.107 & & $\pspar\approx 0.182$ &  0.0681 & & $\pspar\approx 0.501$\\
		\A\PCL & 0.0935 & & & 0.0576 & &\\
		\A\LNL & 0.0438 & & $\scale_1\approx 0.095$ $\scale_5\approx 0.818$ &  0.0225 & &$\scale_1\approx 0.0225$ $\scale_5\approx 1.091$\\
		& & & & & & \\
		\M\MNL & 11.518 & -77.635 & &  5.593 & =0 &\\
		\M\PSL & 12.932 & -91.038 & $\pspar\approx 0.173$ &  8.690 & -47.403 & $\pspar\approx 0.490$\\
		\M\PCL & 12.983 & -108.293  & &  7.809 & -55.479 &\\
		\M\LNL & 0.490 & -142.935 & $\scale_1\approx 73.916$ $\scale_5\approx 1175.47$ & 0.319 & -0.258 &$\scale_1\approx 18.235$ $\scale_5\approx 118.725$\\
		& & & & & & \\
		\MD\MNL & 7.475 & -43.589 & &  4.619 & =0 &\\
		\MD\PSL & 10.613 & -72.889 & $\pspar\approx 0.145$ &  7.381 &  -41.176 &$\pspar\approx 0.432$\\
		\MD\PCL & 14.328 & -126.223 & &  7.396 & -62.199 &\\
		\MD\LNL & 0.475 & -111.01 & $\scale_1\approx 58.367$ $\scale_5\approx 7518.45$ &  0.569 & =0 & $\scale_1\approx 15.265$ $\scale_5\approx 61.098$\\
		\bottomrule
		\end{tabular}
	\label{tab:network}
\end{table}

\begin{figure}[htbp]
\centering
\includegraphics[width=0.49\textwidth]{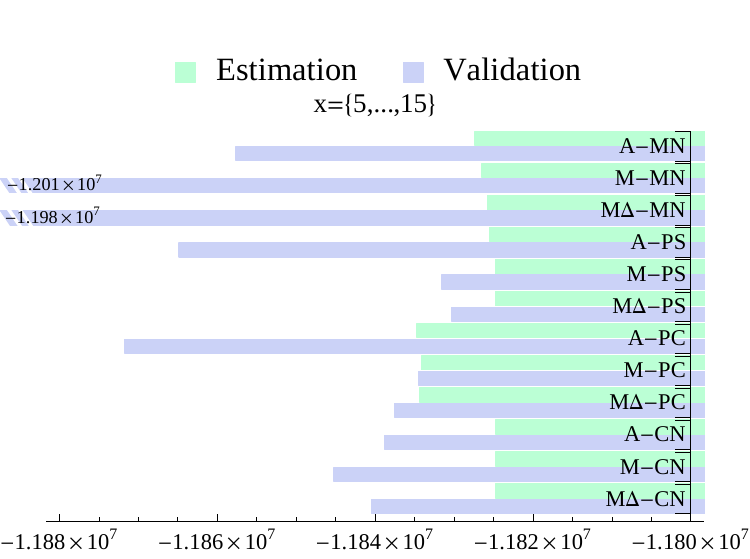}\rule{0.5pt}{5cm}
\includegraphics[width=0.49\textwidth]{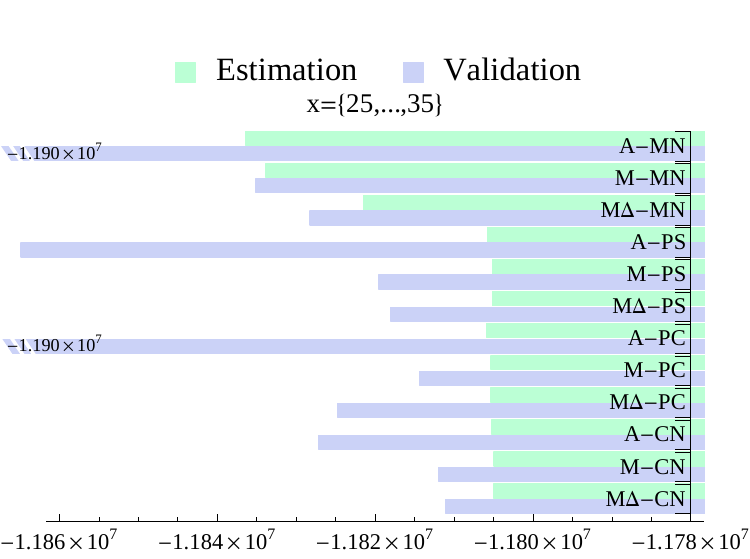}
\caption{Log-likelihoods of model instances for estimation and validation. Parameter estimates of the left estimation are used for the validation in the right, and vice versa.}
\label{fig:LL}
\end{figure}

\begin{landscape}
\begin{figure}[htbp]
\centering
\includegraphics[width=7.5cm]{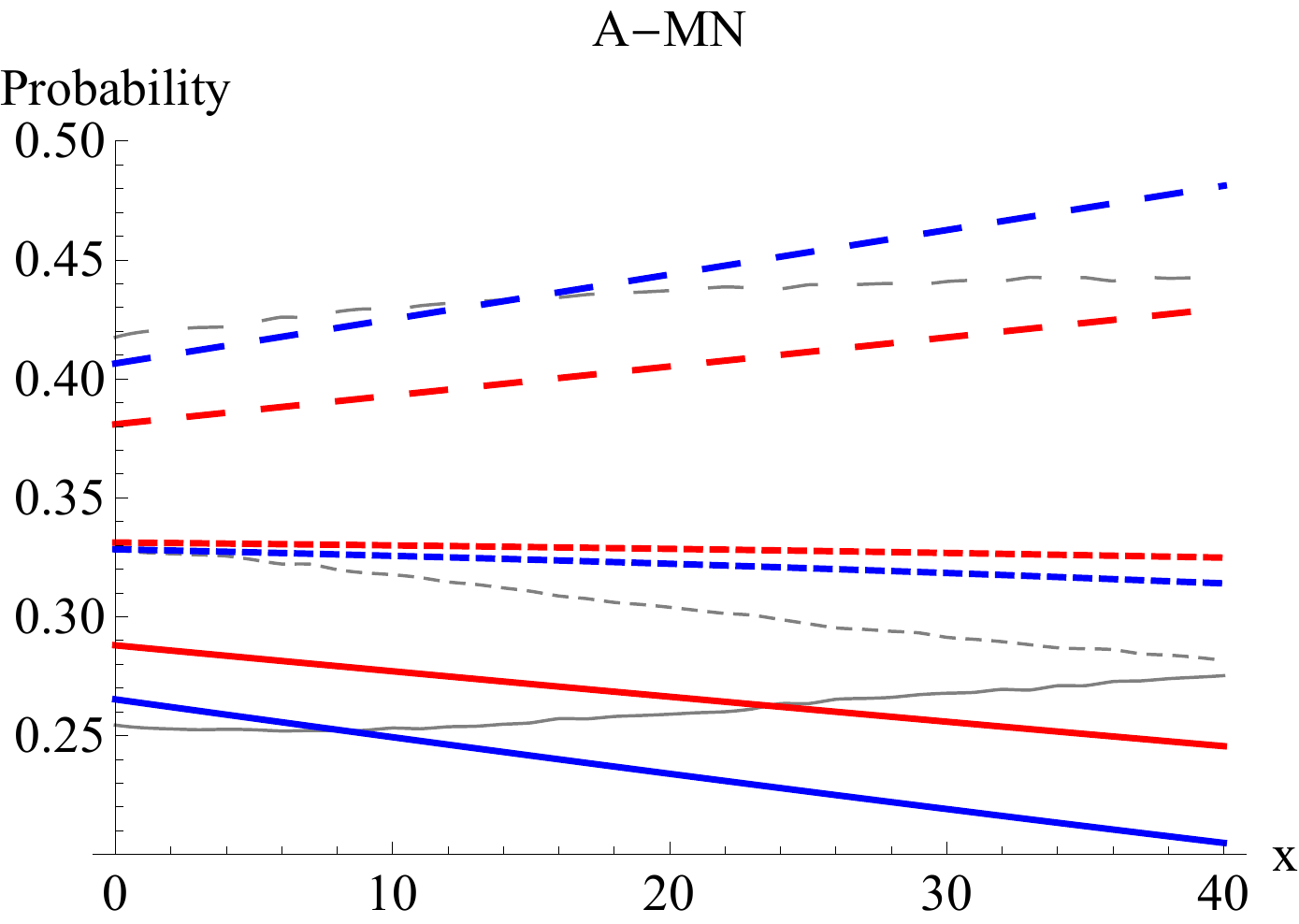}
\includegraphics[width=7.5cm]{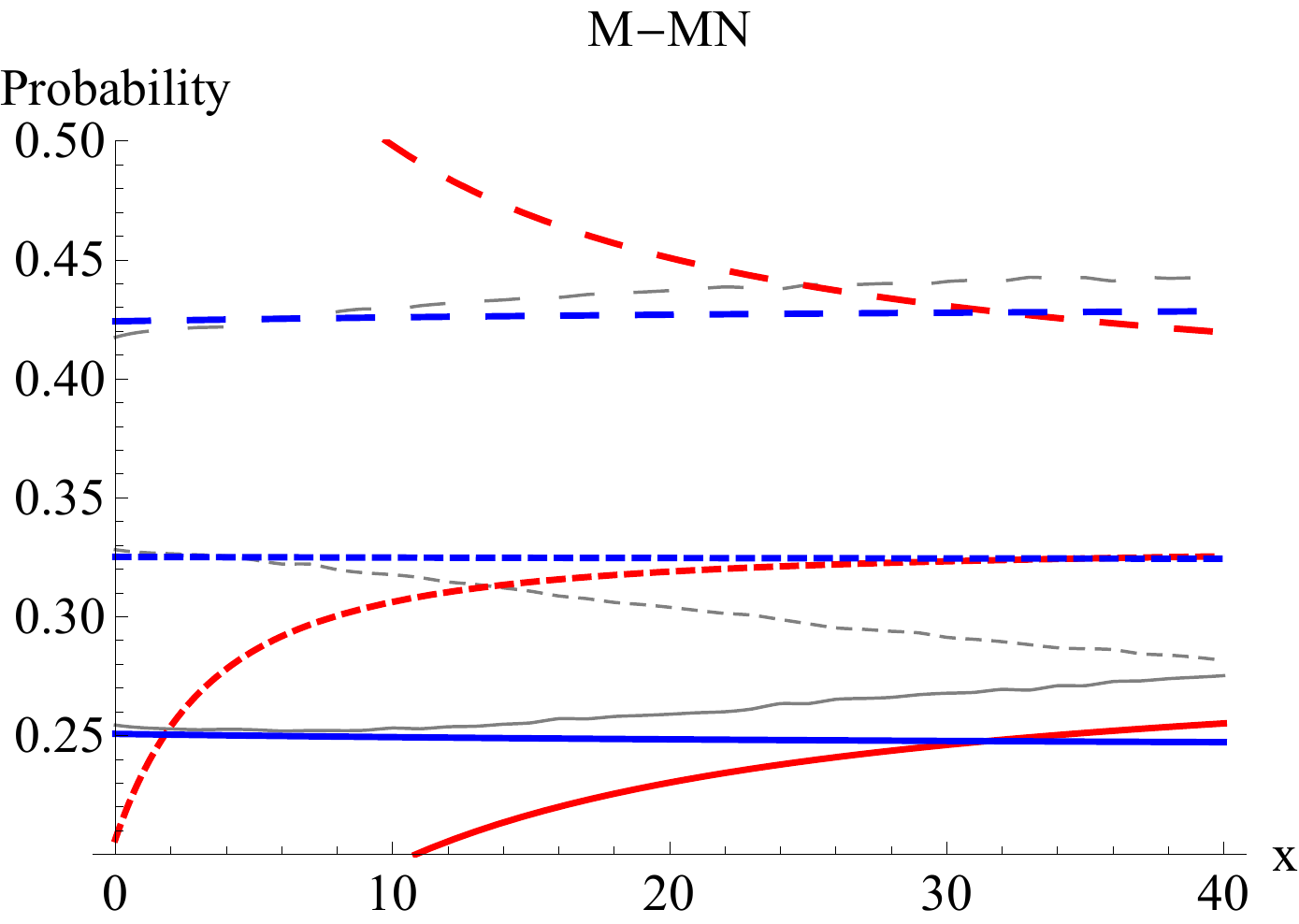}
\includegraphics[width=7.5cm]{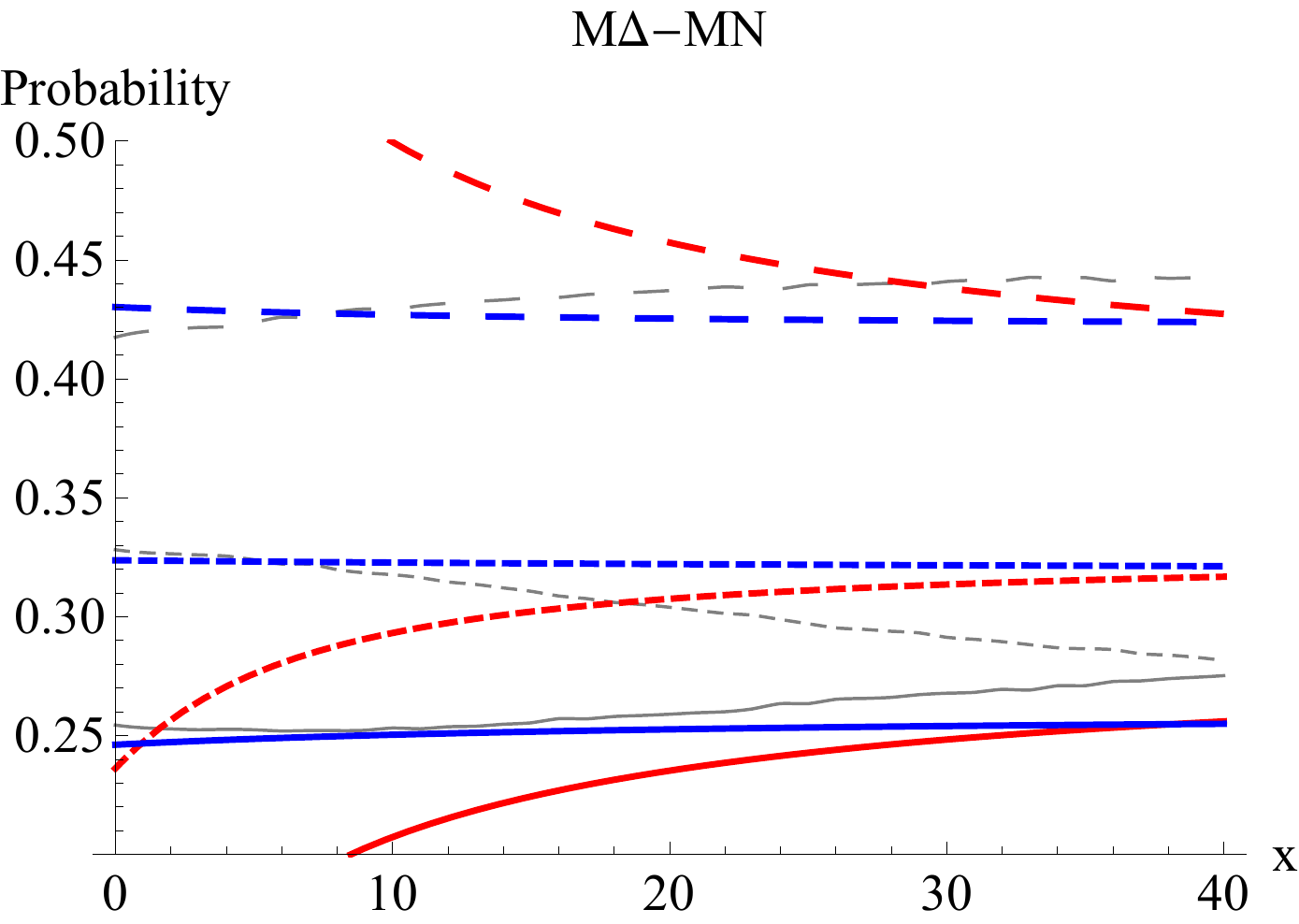}\\[1cm]
\includegraphics[width=7.5cm]{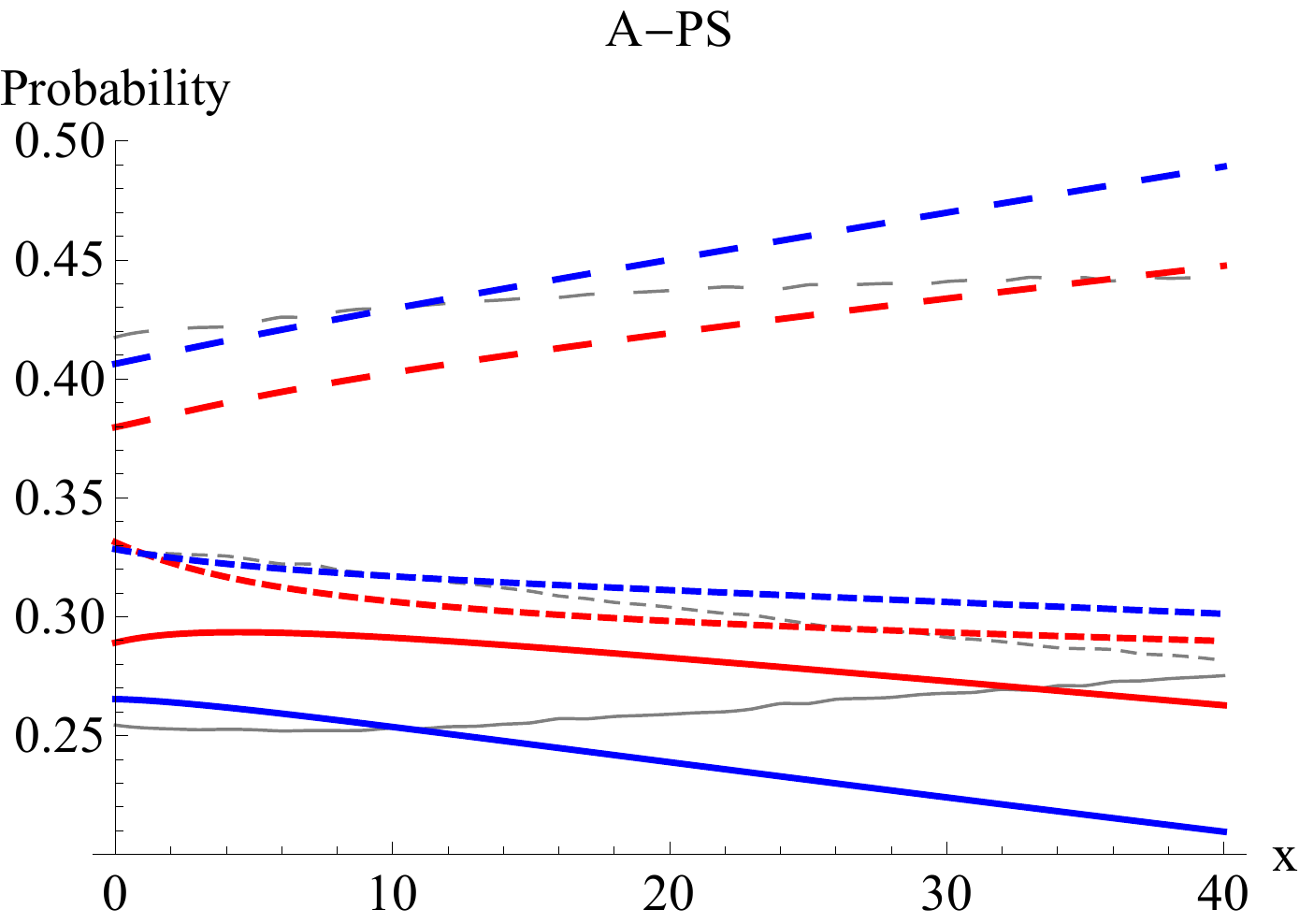}
\includegraphics[width=7.5cm]{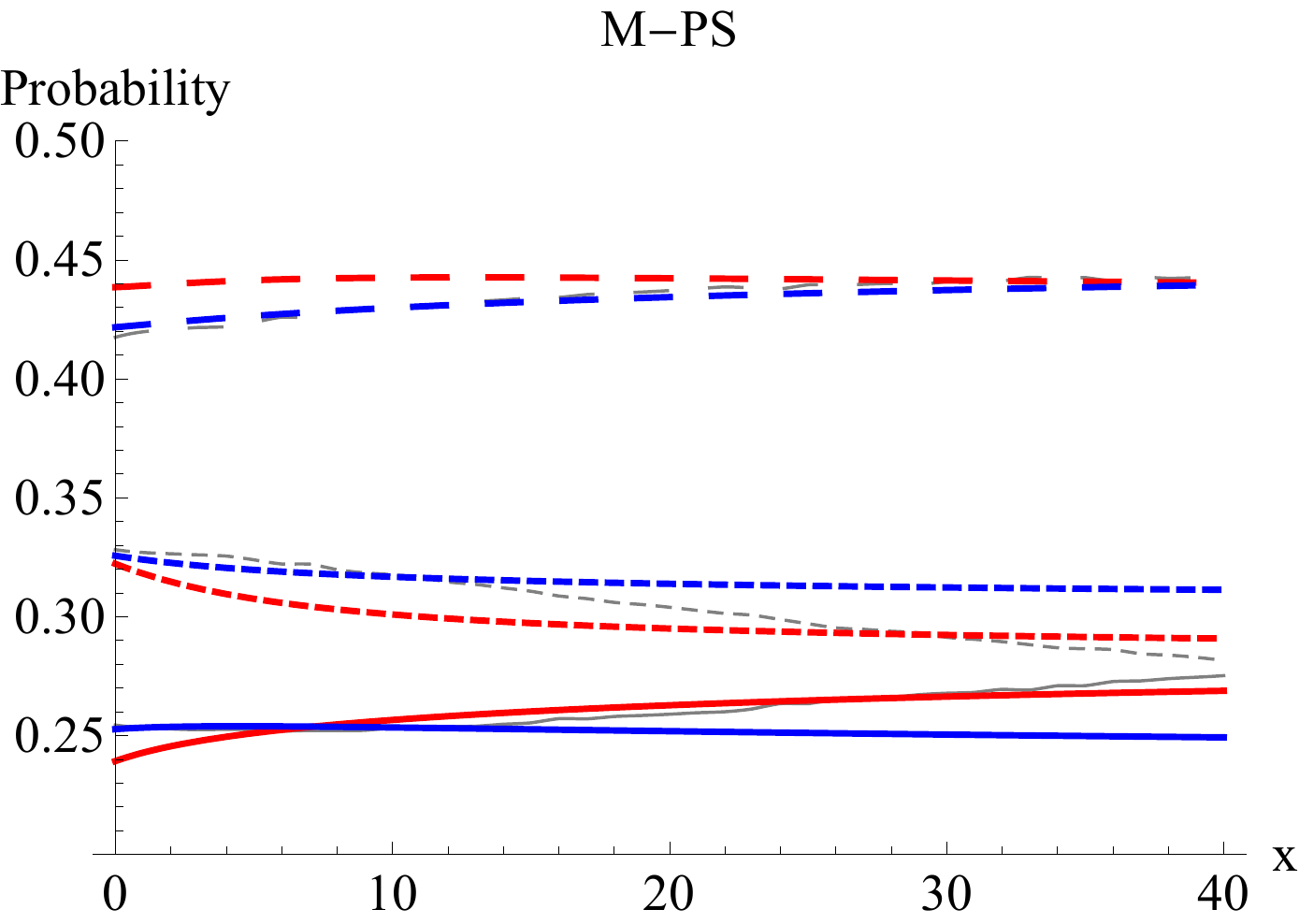}
\includegraphics[width=7.5cm]{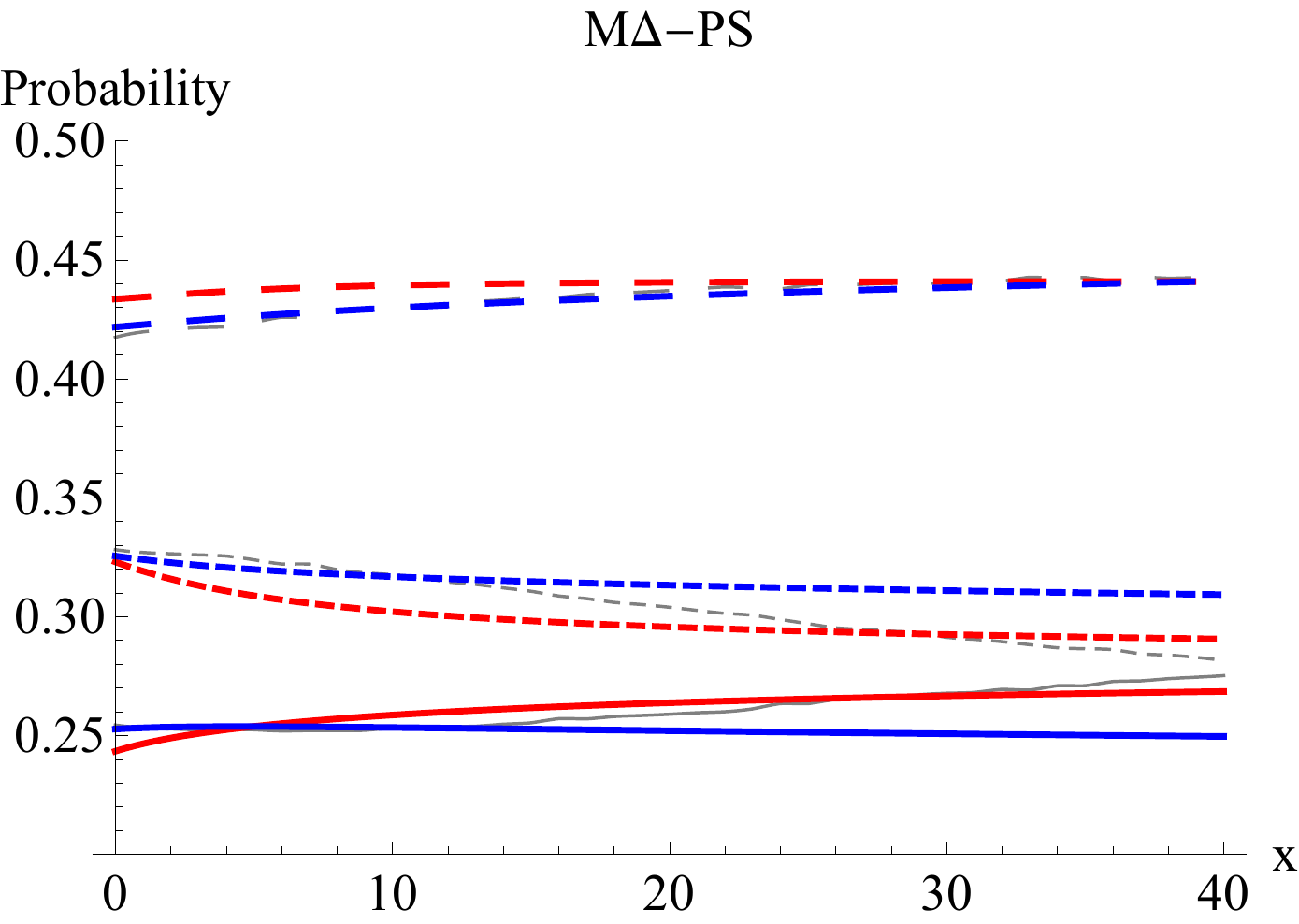}\\[1cm]
\includegraphics[width=22.5cm]{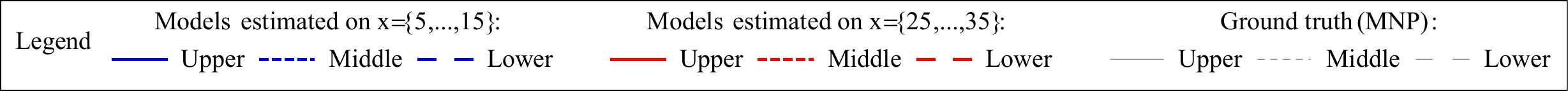}
\caption{Choice probabilities of the \MNL - and \PSL -models.}
\label{fig:MNLPSL}
\end{figure}
\end{landscape}

\begin{landscape}
\begin{figure}[htbp]
\centering
\includegraphics[width=7.5cm]{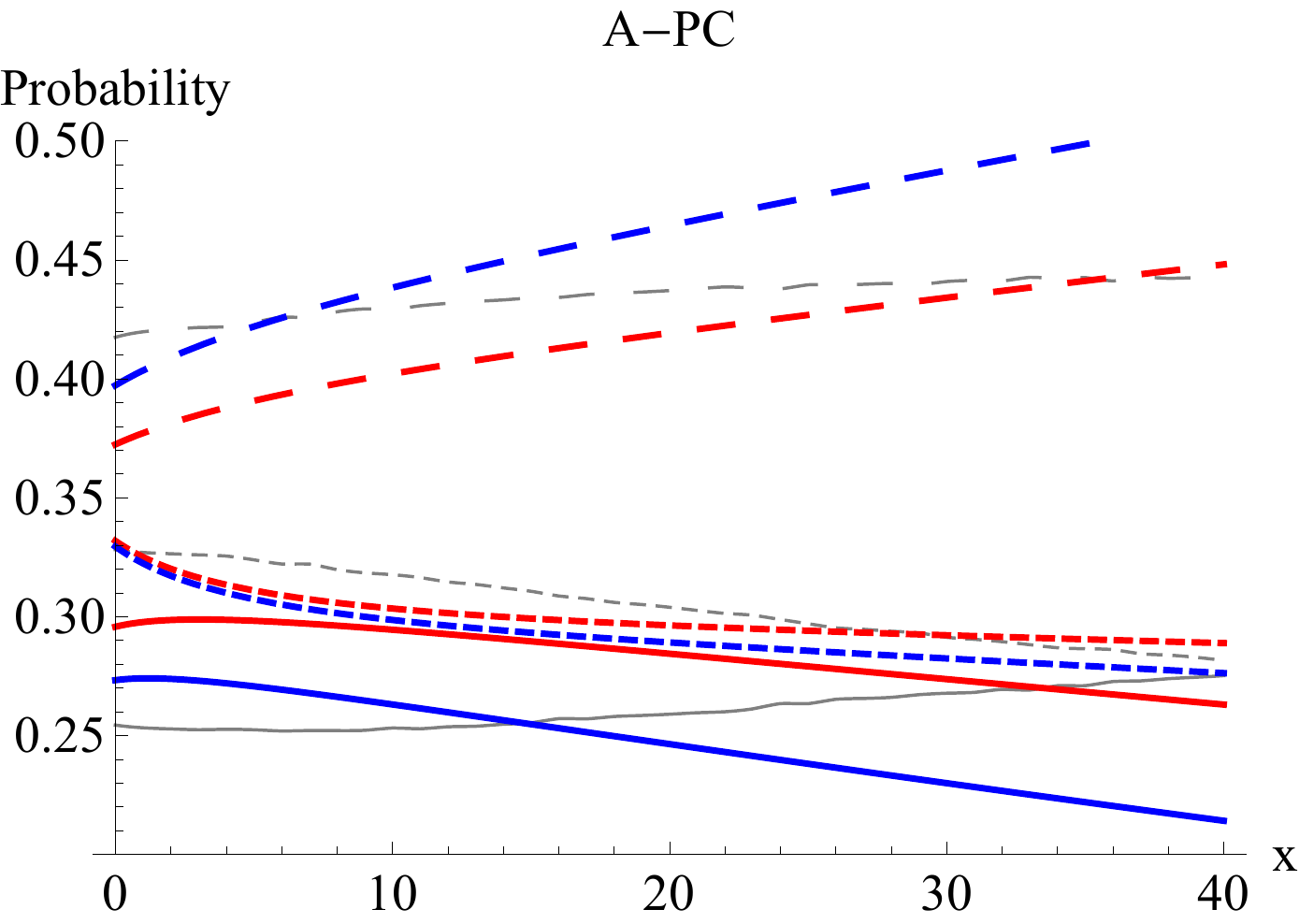}
\includegraphics[width=7.5cm]{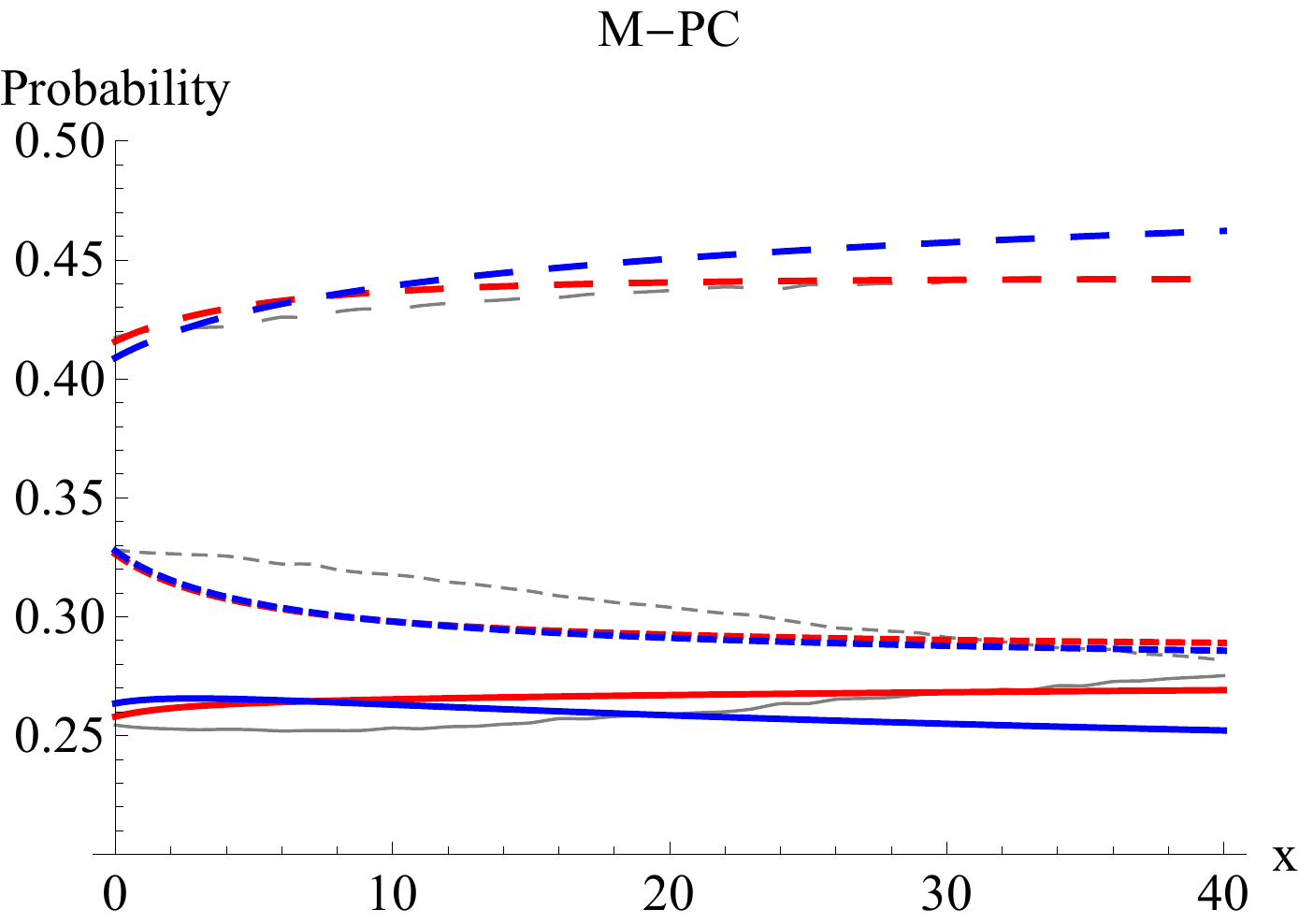}
\includegraphics[width=7.5cm]{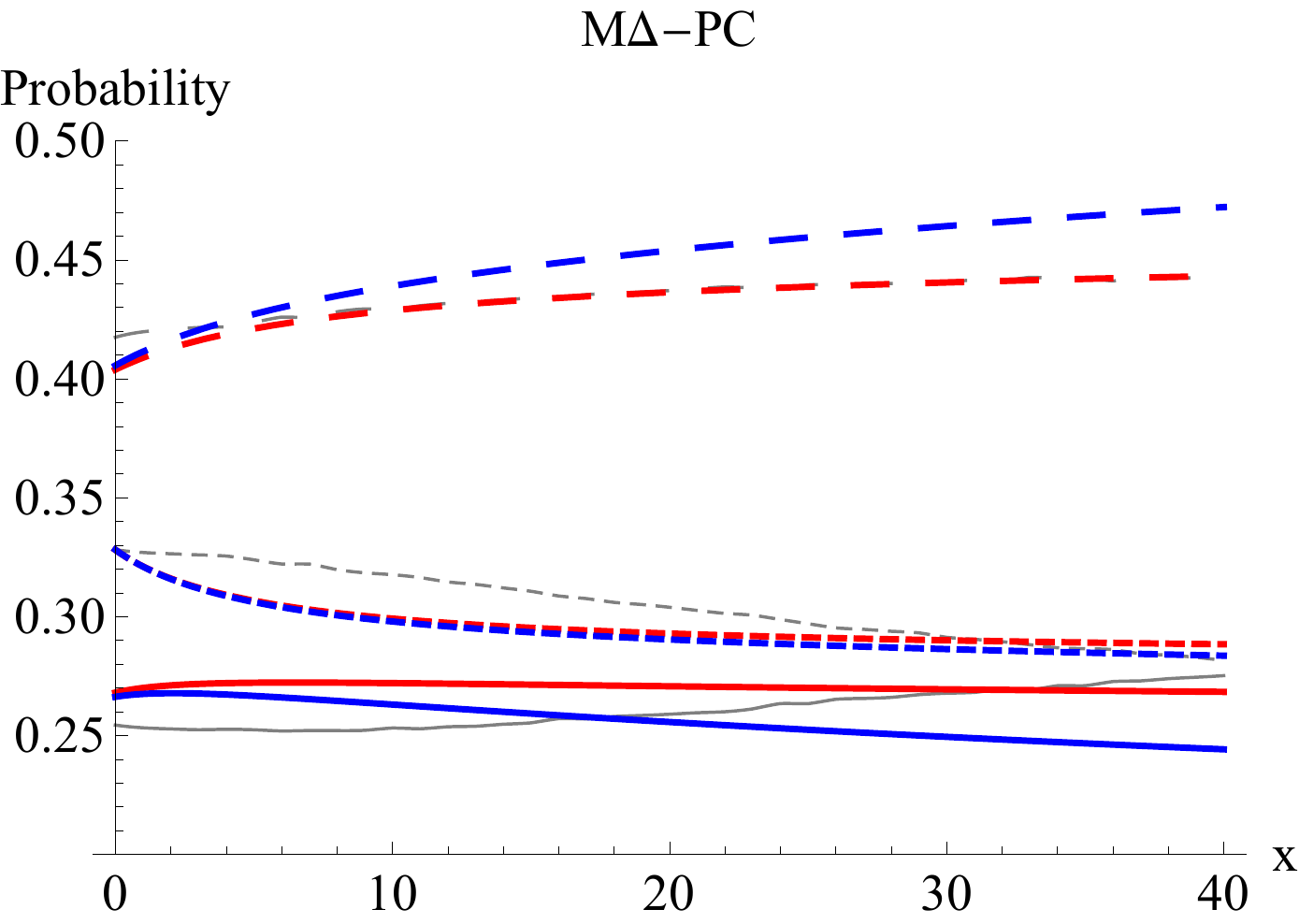}\\[1cm]
\includegraphics[width=7.5cm]{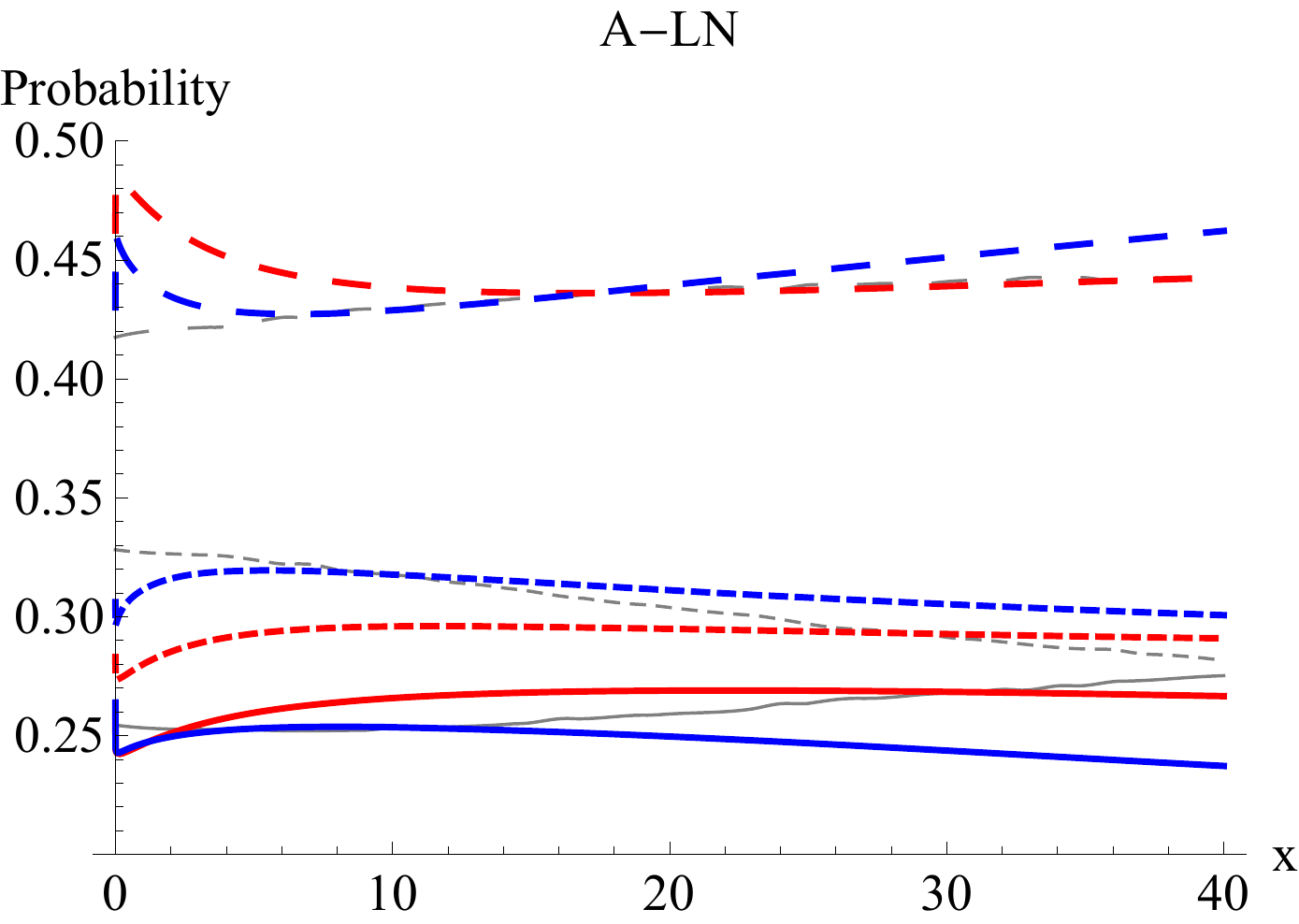}
\includegraphics[width=7.5cm]{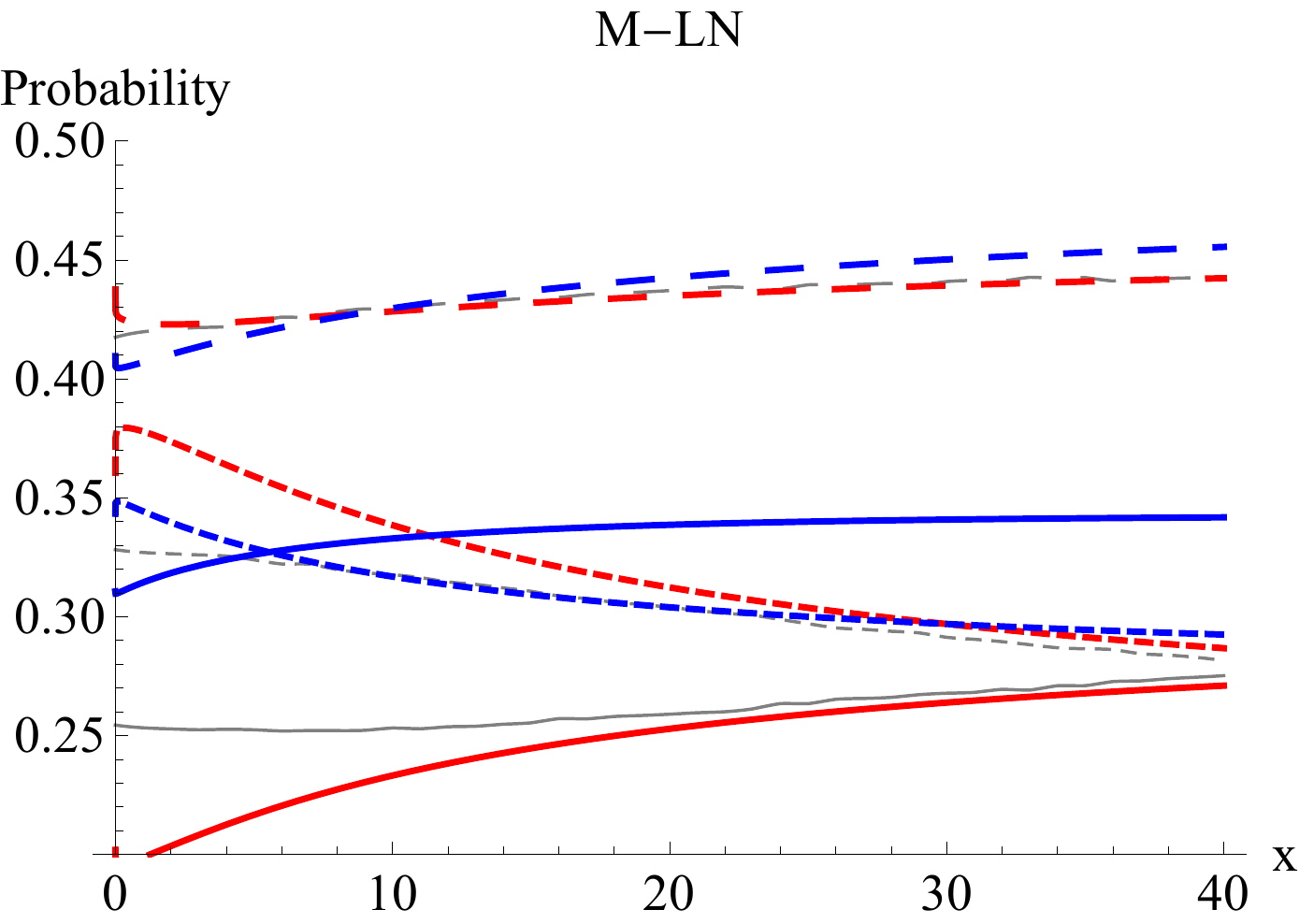}
\includegraphics[width=7.5cm]{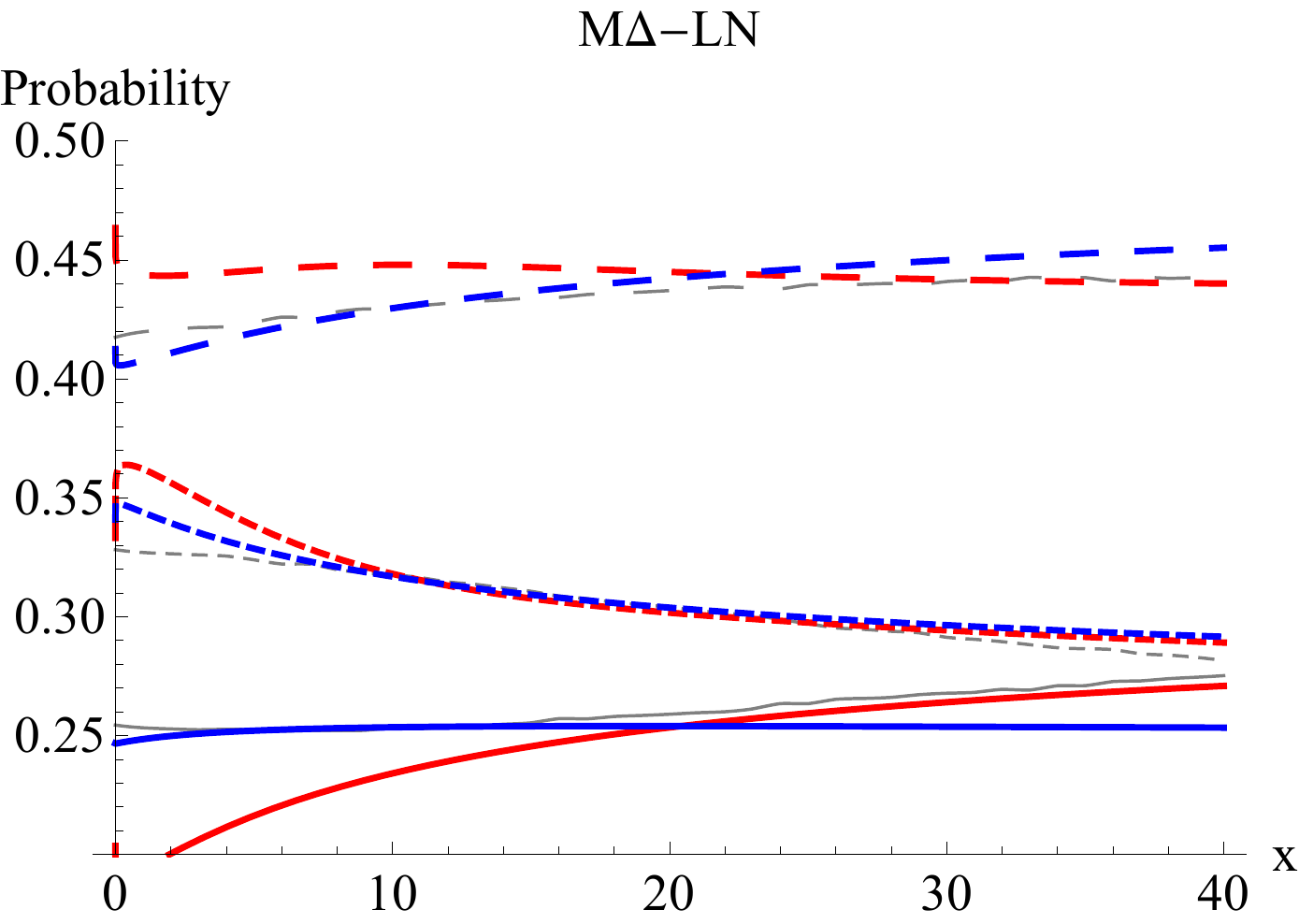}\\[1cm]
\includegraphics[width=22.5cm]{legend.pdf}
\caption{Choice probabilities of the \PCL - and \LNL -models.}
\label{fig:PCLLNL}
\end{figure}
\end{landscape}

Table \ref{tab:network} shows the estimated parameters.\footnote{The q-generalized logit model by \citet{Nakayama2013753} that captures both \A\MNL\ and \M\MNL\ has also been estimated; the results are not presented since the resulting model was always equivalent to the \M\MNL\  model (i.e., not to \A\MNL)} The standard errors are very low due to the large artificial dataset and therefore not reported. Figure \ref{fig:LL} shows the log-likelihoods of the estimation and validation results. Probabilities of the models are found in Figures \ref{fig:MNLPSL} and \ref{fig:PCLLNL}. In every graph two instances of one model and the ground truth are shown, the blue lines depict the route probabilities for the model estimated on $x=\{5,\ldots,15\}$, the red lines depict the route probabilities for the model estimated on $x=\{25,\ldots,35\}$, and the grey lines depict the MNP (i.e., ground truth) probabilities.

\newgeometry{margin=2.5cm}
%\begin{landscape}
%\begin{figure}[htbp]
%\centering
%\includegraphics[width=11cm]{mnw.pdf}\hspace{5mm}
%\includegraphics[width=11cm]{mnwps.pdf}\\[1cm]
%\includegraphics[width=11cm]{qmnl.pdf}\hspace{5mm}
%\caption{Choice probabilities of GEV models. Ground truth probabilities are shown in thin-grey.}
%\label{fig:GEV}
%\end{figure}
%\end{landscape}
\restoregeometry

%To start the analysis, note that \MD\PSL\ is the model with the best estimation and validation results; it is absolutely the best on the $\{5,\ldots,10\}$-data, and has the best validation on the $\{30,\ldots,35\}$-data. On the latter dataset, models \M\LNL, \MD\MNL, and \MD\LNL\ have equally good estimation results\footnote{The first six digits of the log-likelihood are equal.}. Furthermore, the log-likelihood of the \MD\PSL\ validations outperforms many estimation log-likelihoods of other models. Its probability plot shows that both estimated models resemble the ground truth very precisely.
%
%Next, identify that the paired combinatorial and link-nested models  Their validation results are not good in general. The parameters of \M\PCL\ and \MD\PCL\ have very high values (for both the constant and scale) which are difficult to interpret. Despite the very large (negative) constant, the plots show that the non-constant part of the systematic utility does influence the probabilities which is due to the very high scale. The link-nested models estimated on dataset $\{30,\ldots,35\}$ have deviant behaviour for low $x$-values. In any case, it is practically impossible to estimate all link specific nest coefficients for a network.  
%
%Remarkably, the \A\MNL\ model performs well whereas it is simple with only one parameter. Its validation results are the best of all additive models.  

The estimation results show that the \MNL\ models have the poorest log-likelihoods of all on the $\{25,\ldots,35\}$-data; these models cannot address the overlap of routes properly. The validation results of the \MNL\ models are also worse than the other types, so they cannot be transferred between short distance and long distance OD pairs. Another remark is that the \PCL\ models, which are designed to capture overlap, perform poorest on the $\{5,\ldots,15\}$-data; this is due the fact that the covariance between routes is low (because the analyst error is dominating), while the \PCL\ model always imposes dependencies. 

In all cases, the \A models have the worst estimation result compared to the \M\ and \MD models. Remarkably, the multiplicative models also outperform the additive models on the $\{5,\ldots,15\}$-data, where the influence of the foreseen travel time is relatively small. The \MD models have better results on all four \MNL\ and \PSL\ estimations, while the \M models have better results on three out of four \PCL\ and \LNL\ estimations. This is because the additional parameter, the constant $c$, for the multiplicative models, which leads to a better fit.  

The multiplicative errors capture the analyst error with constant $c$ in the systematic utility. This constant is indeed larger for the models estimated on the  $\{5,\ldots,15\}$-data, where the analyst error is dominant. Furthermore, note that all scale parameters $\scale$  are smaller for the $\{25,\ldots,35\}$-data (i.e., the variance is higher); this reflects the heteroscedasticity of route utility.  
 
We would like to mention that the \PCL\ and \LNL\ estimation were problematic on other network configurations that we have tried. For the multiplicative models extremely high scales occurred, which can lead to numerical problems. They also generated very unrealistic probabilities outside the area they were estimated on. As mentioned earlier, it is infeasible to estimate all link-specific nest-scales in large networks for \LNL\ models.   

\section{Conclusions and Discussion}

This paper presented twelve route choice models -- of which seven are new -- in a single framework, and assessed them qualitatively and quantitatively. Choice probabilities for all models have the same closed form expression, namely Equation (\ref{eqn:Gprobs}), based on a generating function and a generating vector. The generating function determines how route overlap is captured, which is either multinomial (\MNL), path-size (\PSL), paired combinatorial (\PCL), or link-nested (\LNL). The generating vector determines the utility formula, which is either additive RUM (\Aa), multiplicative RUM (\Mm), or the newly presented multiplicative RUM based on reference routes that only considers differences between routes (\Md).

For the qualitative assessment a basic structure of utility with random foreseen travel time was presented (see Section \ref{sec:rum}). We base our analysis on two postulates on random travel times of road segments that each lead to a different structure of randomness. Empirical evidence provides the new insight that the foreseen travel time distribution's mean and standard deviation have a linear relationship, contrary to a linear relationship between its mean and variance. The homoscedastic additive models are not able to capture the random foreseen travel time, but multiplicative models do allow for this. Furthermore, differences in normalization, identification, and invariance are pointed out. The constant in systematic utility in multiplicative models does not have to be normalized. This allows more degrees of freedom and a better fit on the data, but makes it more difficult to compare models from the different paradigms directly. One main advantage of the generic GMEV framework is that it can be analysed as a whole, we show this by providing the equivalent stochastic user equilibrium formulation for all models. 

To show the distinctiveness of the \MD models, each model's behaviour under basic network changes was analysed. Only the \MD models can reproduce realistic behaviour when the characteristics of parallel and serial links change.  

To test the models' potential on real networks, and to test whether they can be applied on datasets on which they are not estimated, a carefully constructed network example was presented. Based on our analyses, we expect good performance of the \M\PSL\ and \MD\PSL\ models for route choice on real networks. They can capture overlap sufficiently, and they can handle random foreseen travel time. Also \citet{Fosgerau2009494} and \citet{Chikaraishi2015} have both compared additive with multiplicative formulations on multiple datasets, and they found that the multiplicative models have a better fit for all datasets.\footnote{Note that \citet{Chikaraishi2015} state their model equals logit when $q=0$ and to weibit when $q=1$; however, they do not fully use the flexibility of weibit. When they would have estimated a constant for the weibit models (similar to the approach in this paper), the q-generalized logit models with $0<q<1$ will also become weibit.} However, additional empirical estimation and validation is required to conclusively assess all models. Finally, the \PCL\ and \LNL\ models are problematic to estimate on some other networks we tried, and they can cause numerical problems. 

This paper does not discuss the route generation or sampling related to the explicit route sets of the models. As pointed out earlier, a correct sample of routes is required to obtain unbiased parameter estimates \citep{Frejinger2009984}. For econometrically sound applications, sampling techniques as presented by \citep{Frejinger2009984,Flotterod2013,Guevara2013} have to be adapted for the new models. Models with implicit route sets \citep{Dial197183,Papola2013CACIE,Fosgerau201370,Mai2015} do not have this issue, but they might lead to unrealistic routes. The connection between the link-based MEV model of \citep{Papola2013CACIE} with the route-based GMEV framework does not seem feasible due to the different base units, but if it exists, it might lead to new GMEV model instances. On the other hand, \citet{Prato2012} points out conceptual and empirical reasons that plea for the explicit approach.

%\todosmall{Use less abbr. here, in order create an exective summary in combination with the introduction}
%\todo{OLD}
%The route choice models based on random utility maximization with closed-form expressions for the probabilities are presented and analysed. The heteroscedastic and correlated nature of route utilities is the major challenge in these models. Correlation can be addressed to some degree by the MEV models based on one of these generating functions: $\GPSL$, $\GPCL$, and $\GLNL$. \MNWPS\ is the only GEV model that can capture correlation. heteroscedasticity is resolved in the multiplicative and \MD case MEV models, and also in all GEV models. 
%
%In Theorem \ref{thm:PS} the relation between ordinary models with an additional correction term in the utility formulation and the MEV model \A\PSL\ is provided. It is shown that C-Logit, Path-Size Logit and Path-Size Correction Logit are equivalent to \A\PSL. 
%
%The models \M\PSL, \M\PCL, \M\LNL, \MD\MNL, \MD\PSL, \MD\PCL\ and \MD\LNL\ are presented for the first time. The best performing models on the network example are among these. The additional characteristic of the \MD case MEV models is that the probabilities are determined by the (cost-)ratio between the non-overlapping parts of each route pair. Remember that for multiplicative MEV models the probabilities are determined by the (cost-)ratio between complete routes, and for additive MEV models the probabilities are determined by the difference (in cost) between complete routes. Unfortunately, an additional normalization is required in the \MD case.  

\section{Acknowledgements}
We would like to thank Michel Bierlaire for the discussion with him on an early version of this paper, Kees van Goeverden for preparing the OViN data, and Giselle de Moraes Ramos-Heydendael for providing the foreseen travel time data. The study is performed as part of the Innovative Pricing for Sustainable Mobility (iPriSM) project in the Sustainable Accessibility of the Randstad (SAR) program of the Dutch National Science Foundation (NWO). The main theory in this paper is initiated during a visit of the corresponding author to the Institute of Transport and Logistics Studies at the University of Sydney.  

\bibliographystyle{apalike}
\bibliography{choice}

\end{document}